\input ./style/arxiv-general.cfg
\documentclass[aop,MSNbibl,seceqn,noautosecdot,dvips]{arximspdf}
\makeatletter
   \@ifpackageloaded{graphicx}{}{\usepackage{graphicx}}
\makeatother

\doi{10.1214/15-AOP1053}
\volume{44}
\issue{5}
\pubyear{2016}
\firstpage{3399}
\lastpage{3430}
\docsubty{FLA}

\makeatletter
\def\sfrac#1#2{#1/#2}

\def\sklfrac#1#2{(#1/#2)}

\newcommand{\rrvert}{\vert}
\newcommand{\rrVert}{\Vert}
\newcommand{\llvert}{\vert}
\newcommand{\llVert}{\Vert}
\renewcommand{\mid}{|}
\newcommand{\eqref}[1]{(\ref{#1})}
\newtheorem{theorem}{Theorem}[section]
\newtheorem{lemma}[theorem]{Lemma}
\newtheorem{proposition}[theorem]{Proposition}
\newproclaim{rem}[theorem]{Remark}
\newcommand{\erf}{\operatorname{erf}}
\newcommand{\lip}{\operatorname{lip}}
\newcommand{\ind}{\mathbf{1}}
\newcommand{\grad}{\nabla}
\newcommand{\cG}{{\mathcal G}}
\newcommand{\cA}{{\mathcal A}}
\newcommand{\cE}{{\mathcal E}}
\newcommand{\cH}{{\mathcal H}}
\newcommand{\cC}{{\mathcal C}}
\newcommand{\cT}{{\mathcal T}}
\newcommand{\cB}{{\mathcal B}}
\newcommand{\bsin}{{\overline{\sin}}}
\newcommand{\DD}{\mathrm{d}} 
\newcommand{\bbE}{{\mathbb E}}
\newcommand{\bbP}{{\mathbb P}}
\newcommand{\bbR}{{\mathbb R}}
\newcommand{\bbZ}{{\mathbb Z}}
\newcommand{\gb}{\beta}
\newcommand{\gep}{\varepsilon} 
\newcommand{\gD}{\Delta}
\newcommand{\gO}{\Omega}
\newcommand{\gl}{\lambda}
\makeatother

\begin{document}
\begin{frontmatter}

\title{The cutoff profile for the simple exclusion process on the circle}
\runtitle{Cutoff profile for exclusion on the circle}

\begin{aug}
\author[A]{\fnms{Hubert}~\snm{Lacoin}\corref{}\ead
[label=e1]{lacoin@impa.br}}
\runauthor{H. Lacoin}
\affiliation{IMPA---Instituto Nacional de Matem\'atica Pura e Aplicada}
\address[A]{IMPA---Instituto Nacional de Matem\'atica Pura e Aplicada\\
Estrada Dona Castorina 110\\
22460-320, Rio de Janeiro\\
Brasil\\
\printead{e1}}
\end{aug}

%
\received{\smonth{3} \syear{2015}}
%
\revised{\smonth{7} \syear{2015}}

%
\begin{abstract}
In this paper, we give a very accurate description of the way the
simple exclusion process relaxes to equilibrium.
Let $P_t$ denote the semi-group associated the exclusion on the circle
with $2N$ sites and $N$ particles.
For any initial condition $\chi$, and for any $t\ge\frac{4N^2}{9\pi
^2}\log N$,
we show that
the probability density $P_t(\chi, \cdot)$ is given by an exponential
tilt of the equilibrium measure by the main eigenfunction of the
particle system.
As $\frac{4N^2}{9\pi^2}\log N$ is smaller than the mixing time
which is $\frac{N^2}{2\pi^2}\log N$,
this\vspace*{2pt} allows to give a sharp description of the cutoff profile:
if $d_N(t)$ denote the total-variation distance starting from the worse
initial condition we have
\[
\lim_{N\to\infty} d_N \biggl(\frac{N^2}{2\pi^2}\log N+
\frac
{N^2}{\pi^2}s \biggr)=\operatorname{erf} \biggl(\frac{\sqrt{2}}{\pi
}e^{-s}
\biggr),
\]
where $\operatorname{erf}$ is the Gauss error function.
\end{abstract}

%
\end{frontmatter}

\section{\texorpdfstring{Introduction.}{Introduction}}\label{sec1}

The exclusion process is a lattice interacting particle system where
particles perform independent nearest-neighbor random walks with the
added constraint that each site can be occupied by at most one particle (see the classic references \cite{cfLiggett2} and  \cite{cfLiggett} for a complete introduction to the subject).
It is a very simplified model for a gas of interacting particles. We
consider in this \hyperref[sec1]{Introduction} the case were the lattice
is either a $d$-dimensional torus or hypercube of side length~$N$.
The number of particle is chosen to be proportional to the number of sites.

In this paper, we investigate how the particle system starting
far away from equilibrium, relaxes to
its equilibrium state. This question can in fact be treated with
different point of views:
\begin{itemize}
\item One can describe the evolution of the system at the macroscopic
level, giving the evolution of
the density of particle in space and time. This is the study of
hydrodynamic limits (see \cite{cfKL} for an introduction to the subject).
\item One can adopt a microscopic point of view, and look at the
evolution of the law of the particle, and in particular,
its total variation distance to the equilibrium law. This is the study
of the Markov chain's mixing time (see \cite{cfLPW}).
\end{itemize}
With the first point of view, it is now well known that under diffusive
rescaling (space rescaled by $N$ and
time rescaled buy $N^2$), the density of particle evolves like the
solution of the heat equation.
The result is valid in any dimension (see \cite{cfKL} for references).

Concerning the mixing-time approach, progresses are more recent. It has
been shown by Morris that in any dimension the time needed to come
close to equilibrium in total variation was of order $N^2\log N$ \cite{cfMor}.
In dimension $1$, more refined estimates have been obtained and gave
the exact location of the mixing time either for the segment
\cite{cfLac1} or the circle \cite{cfLac2} with lower bounds proved
earlier by Wilson \cite{cfWilson} (see also \cite{cfRO} for
results in the case of arbitrary graph, and \cite{cfDS2,cfLL}).
A natural question is then what does the law of the particle system
look like when it is about to reach equilibrium.

At equilibrium, the law of the distribution is uniform over all
particle configurations.
Another way to see it is to say that the state of each site (occupied
or vacant) is given by
a field of i.i.d. Bernoulli variables whose sum is conditioned to be
equal to the number of particle (which is a fixed parameter).

What would be natural to expect then, is that up to a small correction,
before equilibrium, the particle distribution still is a
conditioned product measure, but that the Bernoulli variables are no
more identically distributed: there is a space dependent bias
which is given by the solution of the heat equation. This brings a
strong connection between the problem of the mixing time and that
of the
hydrodynamical limit. This connection was previously underlined by Lee
and Yau when studying the related issue of $\log$-Sobolev constant for
the simple exclusion
\cite{cfLY}.
Indeed in the case of small bias, with some minor efforts one can
derive sharp estimates on the total-variation distance between
the conditioned product of biased Bernoulli and the equilibrium
measure. This can be turned into a precise prediction on
how the total-variation distance drops from one to zero, the cutoff profile.
The present paper brings this heuristic picture on a rigorous ground in
the case of the exclusion on the circle.

\section{\texorpdfstring{Model and results.}{Model and results}}

\subsection{\texorpdfstring{The process.}{The process}}

We consider $\bbZ_{2N}:=\bbZ/(2N\bbZ)$, the discrete circle with
$2N$ sites and we place $N$ particles
on it, with \textit{at most} one particle per site. With a slight
abuse of notation,
we will sometimes use elements of $\{1,\dots,2N\} \subset\bbZ$ to
refer to elements of $\bbZ_{2N}$.

The exclusion process on $\bbZ_{2N}$ is a dynamical evolution of the
particle system which can be described informally as follows:
each particle tries to jump independently on its neighbors with
transition rates $p(x,x+1)=p(x,x-1)=1$,
but the jumps are canceled if
a particle tries to jump on a site which is already occupied.

Let us describe the chain more formally.
We adopt the convention that $1$ denotes a particle and $-$1 denotes an
empty site. This is not the most usual one (empty sites are more often
denoted by $0$)
but it proves to be more practical in our computations.
Our state-space is defined by
%
\begin{equation}
\gO_{N}= \Biggl\{ \eta\in\{-1,1\}^{\bbZ_{2N}} \Big| \sum
_{x=1}^{2N} \eta(x)= 0 \Biggr\}.
\end{equation}
Given $\eta\in\gO$ define $\eta^x$ the configuration obtained by
exchanging the content of site $x$ and $x+1$
%
\begin{equation}
\cases{ \eta^x(x):=\eta(x+1),
\cr
\eta^x(x+1):=\eta(x),
\cr
\eta^x(y)=\eta(y)\qquad\forall y\notin\{x,x+1\}. }
\end{equation}
The exclusion process on $\bbZ_{2N}$ with $N$ particle is the
continuous time Markov process on
$\gO_{N}$ whose generator is given by
%
\begin{equation}
\label{crading} (\mathcal L_N f) (\eta):=\sum
_{x\in\bbZ_{2N}} f\bigl(\eta^x\bigr)-f(\eta).
\end{equation}
The chain is irreducible and reversible, and
the unique invariant probability measure is the uniform probability
measure on $\gO_N$ which we denote by $\mu_N$.
Given $\chi\in\gO_N$ we let $(\eta^{\chi}_t)_{t\ge0}$ denote the
trajectory of the Markov chain starting from~$\chi$. We write
$\bbP [ (\eta^{\chi}_t)_{t\ge0}\in\cdot ]$ for the law
of $(\eta^{\chi}_t)_{t\ge0}$.
We\vspace*{1pt} let $P_t$ denote the Markov semi-group and write $P^\chi_t$ for the
probability measure $P_t(\chi,\cdot)$, $\chi\in\gO_N$.

We measure the distance to equilibrium in terms of total variation distance.
If $\alpha$ and $\beta$ are two probability measures on $\gO$, the
total variation distance between
$\alpha$ and $\gb$ is defined to be
%
\begin{equation}
\label{tv} \llVert \alpha-\beta\rrVert _{\mathrm{TV}}:=\frac{1}{2}\sum
_{\omega\in\gO} \bigl\llvert \alpha (\omega)-\beta(\omega)
\bigr\rrvert =\sum_{\omega\in\gO} \bigl(\alpha(\omega )-\beta(
\omega)\bigr)_+,
\end{equation}
where $x_+=\max(x,0)$ is the positive part of $x$. It measures how
well one can couple two variables with law $\alpha$ and $\beta$.
We define the distance to equilibrium of the Markov chain to be
%
\begin{equation}
\label{dat} d^{N}(t):=\max_{\chi\in\gO_{N}} \bigl\llVert
P^\chi_t-\mu\bigr\rrVert _{\mathrm{TV}}.
\end{equation}

In a previous paper \cite{cfLac2}, we described in detail the
asymptotic behavior of $d^N(t)$.
We proved that around a time of order $\frac{N^2}{2\pi^2}\log N$
the distance to equilibrium drops from $1$ to $0$ in a time window of
width $N^2$,
%
\begin{eqnarray}
\lim_{s\to\infty}\limsup_{N\to\infty}d^{N}
\biggl(\frac{N^2}{2\pi
^2}\log N+ N^2s \biggr)&=&0,
\nonumber\\[-8pt]\\[-8pt]\nonumber
\lim_{s\to-\infty}\liminf_{N\to\infty} d^{N}
\biggl(\frac
{N^2}{2\pi^2}\log N- N^2s \biggr)&=&1.
\end{eqnarray}

The aim of this paper is to complete this picture by identifying, in an
acute way, the pattern of relaxation to equilibrium.
In particular, we are interested in proving the existence and finding
an expression for limiting profile
%
\begin{equation}
\lim_{N\to\infty} d^{N} \biggl(\frac{N^2}{2\pi^2}\log N+
N^2s \biggr).
\end{equation}

To reach this goal, we have to understand what the distribution $P^\chi
_t$ looks like much before the time $\frac{N^2}{2\pi^2}\log N$.

\subsection{\texorpdfstring{The mixing time profile.}{The mixing time profile}}

The main achievement of our paper is to determine the cutoff profile.

\begin{theorem}\label{maintheorem}
The total-variation distance to equilibrium from the worst initial
condition has the following asymptotic profile:
for any $s\in\bbR$, we have
%
\begin{equation}
\lim_{N\to\infty} d^{N} \biggl(\frac{N^2}{2\pi^2}\log N +
\frac
{N^2}{\pi^2} s \biggr)=\erf \biggl(\frac{\sqrt{2}}{\pi}e^{-s} \biggr),
\end{equation}
where $\erf$ is the Gauss error function
%
\begin{equation}
\erf(t):=\frac{2}{\sqrt{\pi}} \int_0^t
e^{-u^2}\,\DD u.
\end{equation}
\end{theorem}

The method by which we obtain the result gives us in fact much more
information about the relaxation of the system:
we are able to characterize fully how $P_t^\chi$ looks like much
before the mixing time, for all initial condition $\chi\in\gO_N$.

\begin{rem}
The fact that the profile of the cutoff is given by a function of the
type $\erf (Ae^{-s} )$ (the constant is not essential since
it depends on the particular process and
the choice for the normalization)
is given by Wilson \cite{cfWilson} as a conjecture
(supported by numerical evidences) for a process very much related to
the exclusion: the adjacent transposition shuffle.
The reason why the function $\erf$ appears is that the last statistic
that comes to equilibrium for the process
(here the first Fourier coefficient of $\eta$, see below)
is well approximated by a Gaussian; the exponential terms are present
because the mean of this Gaussian converges exponentially to zero.
This is a property which is believed to be shared by many Markov chains
and rigorously
known, for example, the random walk on the hypercube \cite{cfDGM}.
Let us mention that however there are the Markov chains which exhibit
cutoff and do not have this property. This is, for instance,
the case of top to random shuffle \cite{cfDFP}, and also of the
transposition shuffle for which
the important statistic, the number of fixed point, behaves like a
Poisson variable (see, e.g., \cite{cfMatthews}).
\end{rem}

\subsection{\texorpdfstring{The description of $P_t^\chi$ much before equilibrium.}{The description of Ptchi much before equilibrium}}\label{intro}

The main result of the paper, from which we deduce Theorem~\ref
{maintheorem} requires some notation to be introduced.
The time evolution of the density of particles is given by the discrete
heat equation on $\bbZ_{2N}$ and for this reason,
the eigenfunction of the discrete Laplacian on the circle plays an
important role in our analysis; in particular, those in the eigenspace
corresponding to the spectral gap:
$x\mapsto\cos (\frac{\pi x}{N}  )$, and $x\mapsto\sin
 (\frac{\pi x}{N}  )$.

To describe the projection of $\chi\in\gO_N$ on this eigenspace, it
is more convenient for us to have one positive coefficient than two
real ones, and for this reason we introduce
$\theta(\chi)$ which is the ``phase'' of $\chi$ in the first
eigenspace. It is the unique $\theta\in[0,2\pi)$ satisfying
%
\begin{eqnarray}
\label{system} \sum_{x\in\bbZ_{2N}} \chi(x) \cos \biggl(
\frac{\pi x}{N}+\theta \biggr) &=&0,
\nonumber\\[-8pt]\\[-8pt]\nonumber
\sum_{x\in\bbZ_{2N}} \chi(x) \sin \biggl(\frac{\pi x}{N}+
\theta \biggr) &>& 0,
\end{eqnarray}
or $\theta(\chi)=0$ if the system has no solution. We denote by
$b(\chi)$ the first Fourier coefficient
of $\chi$
%
\begin{equation}
b(\chi):=\frac{1}{N}\sum_{z\in\bbZ_{2N}} \chi(x) \sin
\biggl(\frac{\pi x}{N}+\theta \biggr).
\end{equation}
Note that $b(\chi)=0$ in the case where \eqref{system} has no
solution. In the case where $\chi=+1$ for $x\in\{1,\dots,N\}$ and
$-1$ elsewhere,
$\eta(\chi)=\frac{\pi}{2N}$.

If $\mu$ is a probability measure on a state-space $\gO$ and that $f$
is a function $\gO\to\bbR$,
we use the following notation for the expectation:
%
\begin{equation}
\nu(f):=\nu\bigl(f(\eta)\bigr):=\sum_{\eta\in\gO} f(\eta)
\nu(\eta).
\end{equation}
We define given $N$, $\alpha>0$ and $\theta\in[0,2\pi)$. We define
$\nu^{N,\alpha,\theta}$ to be the probability measure on
$\gO_N$ with the following Radon--Nikodym density:
%
\begin{equation}
\frac{\DD\nu^{N,\alpha,\theta}}{\DD\mu_N}(\eta):= \frac{e^{\alpha a_\theta(\eta)}}{\mu_N (e^{\alpha a_\theta
(\eta)} )},
\end{equation}
where
%
\begin{equation}
a_\theta(\eta):=\sum_{x\in\bbZ}\eta(x) \sin
\biggl(\frac{\pi
x}{N}+\theta \biggr).
\end{equation}
Finally, let us set
%
\begin{equation}
\gl_N:=2 \biggl(1-\cos \biggl( \frac{\pi}{N} \biggr) \biggr).
\end{equation}
We note that $\gl_N$ is the spectral gap of the simple random walk on
$\bbZ_{2N}$ (with jump rate one in each direction),
and hence from \cite{cfCLR}, Section~4.1.1, it is also the spectral
gap of the exclusion process on $\bbZ_{2N}$.

The main result of the paper tells us that much before the mixing time,
$P^\chi_t$
is close to $\nu^{N,\alpha,\theta}$ for an appropriate choice of
$\alpha$ and $\theta$.

\begin{proposition}\label{mainresult}
For all $N$ sufficiently large, for all $\chi\in\gO_N$
for all $t\ge\frac{4\pi^2}{9N^2}$, we have
%
\begin{equation}
\bigl\llVert P^\chi_t-\nu^{N,b(\chi)e^{-\gl_N t}, \theta(\chi)} \bigr\rrVert
_{\mathrm{TV}}\le(\log\log N)^{-1}.
\end{equation}
\end{proposition}

Theorem~\ref{maintheorem} follows from Proposition~\ref{mainresult}
by using the following lemma.

\begin{lemma}\label{tvnu}
For all $K>0$,
for all $N$ sufficiently large (depending on $K$),
%
\begin{equation}
\lim_{N\to\infty} \mathop{\sup_{\gamma\in[0,K]}}_{\theta\in
[0,2\pi)}
\biggl\llvert \bigl\llVert \nu^{N,\gamma N^{-1/2}, \theta}-\mu_N \bigr\rrVert
_{\mathrm{TV}}-\erf \biggl(\frac{\gamma}{\sqrt{8}} \biggr) \biggr\rrvert =0.
\end{equation}
\end{lemma}

\subsection{\texorpdfstring{Exclusion with an arbitrary number of particle.}{Exclusion with an arbitrary number of particle}}

We have chosen to present here the result only in the case where the
number of particles is equal to half of the number of sites.
However, \textit{mutatis mutandis}, the proof adapts to the case of
$k_N$ particle $k_N\le N$ on $\bbZ_{2N}$ where $k_N$ tends to infinity
with $N$ (the case $k\ge N$ can be treated by symmetry). Let us discuss
here what the results are in that case and how they can be obtained.

While the case of $k_N$ proportional to $N$ can be derived directly
from the proof presented here, it turns out that
some of the technical lemmas (e.g., Proposition~\ref{iamzea}) breaks
down if $k_N$ grows much slower, that is, like $\log N$.
However, adapting the techniques developed specifically for the case of
slowly growing $k_N$ in \cite{cfLac2},
all technical obstacles
can be overcome.

To close this discussion, let us mention what the cutoff profiles are
in those cases.
When $k_N=\lceil\alpha N \rceil$ for some $\alpha=(0,1)$, we have
%
\begin{equation}
\lim_{N\to\infty} d^{N} \biggl(\frac{N^2}{2\pi^2}\log N +
\frac
{N^2}{\pi^2} s \biggr)=\erf \biggl(\frac{\sin (\sfrac{\alpha
\pi}{2} )}{\pi\sqrt{\alpha(1-(\alpha/2))}}e^{-s}
\biggr).
\end{equation}
When $k_N$ satisfies $1\ll k_N\ll N$, we have
%
\begin{equation}
\lim_{N\to\infty} d^{N} \biggl(\frac{N^2}{2\pi^2}\log
k_N + \frac
{N^2}{\pi^2} s \biggr)=\erf \biggl(\frac{1}{2}e^{-s}
\biggr).
\end{equation}

\subsection{\texorpdfstring{Organization of the paper.}{Organization of the paper}}
In Section~\ref{profile}, we prove Theorem~\ref{maintheorem} from
Proposition~\ref{mainresult}, and also give a
proof of Lemma~\ref{tvnu}.
In Section~\ref{decompo}, we decompose the proof of Proposition~\ref
{mainresult} into three key statements,
whose proofs are, respectively, given in Sections~\ref{fluctutu},~\ref
{coupel} and~\ref{relax}.

\section{\texorpdfstring{The cutoff profile.}{The cutoff profile}}\label{profile}

\subsection{\texorpdfstring{Proof of Theorem \protect\ref{maintheorem}.}{Proof of Theorem 2.1}}

Let $s\in\bbR$ be fixed. It is straightforward to check that
for
\[
t_{s,N}:= \frac{N^2}{2\pi^2}\log N+ \frac{N^2}{\pi^2} s %
\]
we have
%
\begin{equation}
\sup_{N} \sup_{\chi\in\gO_N} \sqrt{N} b(
\chi)e^{-\gl_N
t_{s,N}}<\infty.
\end{equation}
Hence, using the triangular inequalities, Proposition~\ref{mainresult}
and Lemma~\ref{tvnu} we have
for all $\chi\in\gO_N$,
%
\begin{equation}
\lim\biggl\llvert \bigl\llVert P_{t_{s,N}}^{\chi}-
\mu_N\bigr\rrVert _{\mathrm{TV}}-\erf \biggl(\frac{b(\chi)\sqrt{N}e^{-\gl_N t_{s,N}}}{\sqrt{8}}
\biggr) \biggr\rrvert =0.
\end{equation}
The asymptotic for $d_N(t_{s,N})$ follows if one can identify $\chi$
which maximizes $b(\chi)$.
A few seconds of thoughts are enough to realize that $\chi_{\max}$
defined as follows
is the unique maximizer up to translation:
%
\begin{equation}
\chi_{\max}(x)= \cases{ +1, &\quad for $x=1,\dots, N$,
\cr
-1, &\quad
for $x=N+1,\dots, 2N$.}
\end{equation}
The asymptotic behavior of $b(\chi_{\max})$ is given by the following limit:
%
\begin{equation}
\lim_{N \to\infty} \frac{1}{N}\sum
_{x\in\bbZ_{2N}} \chi_{\max
}(x)\sin \biggl( \frac{x\pi}{N}-
\frac{\pi}{2N} \biggr)=\frac
{4}{\pi}.
\end{equation}
As for any $s\in\bbR$, we have also
%
\begin{equation}
\lim_{N \to\infty} \sqrt{N}e^{-\gl_N t_{s,N}}=e^{-s}
\end{equation}
the result follows from
the continuity of the error function.


\subsection{\texorpdfstring{Proof of Lemma \protect\ref{tvnu}.}{Proof of Lemma 2.4}}
The underlying idea is quite simple: we want to prove that
asymptotically under $\mu_N$, once rescaled
\[
a_{\theta}(\eta):=\sum_{x\in\bbZ}\eta(x) \sin
\biggl(\frac{\pi
x}{N}+\theta \biggr),
\]
converges to a Gaussian.

\begin{lemma}\label{GTLC}
The following statements hold true:
\begin{longlist}[(ii)]
\item[(i)] For a fixed $\theta\in[0,2\pi)$. The quantity
$N^{-1/2}a_{\theta}(\eta)$ converges in law to a standard Gaussian.
Moreover, the convergence is uniform in $\theta$, in the sense that\vadjust{\goodbreak}
for any bounded continuous function $F$ the convergence
%
\begin{equation}
\lim_{N\to\infty} \mu_N \biggl[ F \biggl(
\frac{a_{\theta}(\eta
)}{\sqrt{N}} \biggr) \biggr]= \frac{1}{\sqrt{2\pi}}\int F(u)e^{-\sfrac{u^2}{2}}\,\DD
u,
\end{equation}
holds uniformly in $\theta$.
\item[(ii)]
Moreover, $a_{\theta}(\eta)$ is exponentially concentrated in the
sense that there exists a constant $c>0$ such that
%
\begin{equation}
\mu_N \bigl( \bigl\llvert a_{\theta}(\eta)\bigr\rrvert \ge
\sqrt{N}u \bigr)\le2e^{-cu^2}.
\end{equation}
\end{longlist}
\end{lemma}

Let us explain how we deduce Lemma~\ref{tvnu} from these facts.
We note that
%
\begin{equation}
\bigl\llVert \nu^{N,\gamma N^{-1/2},\theta}-\mu_N\bigr\rrVert _{\mathrm{TV}} =
\frac{1}{2}\mu_N \biggl( \biggl\llvert \frac{e^{\gamma N^{-1/2}a_\theta
(\eta)}}{\mu_N (e^{\gamma N^{-1/2}a_\theta(\eta)} )}-1
\biggr\rrvert \biggr).
\end{equation}
Because of convergence in probability and exponential tightness,
we have
%
\begin{equation}
\lim_{N\to\infty} \mu_N \bigl(e^{\gamma N^{-1/2}a_\theta(\eta
)} \bigr)=
e^{\gamma^2/2}.
\end{equation}
Thus, $\llVert  \nu^{N,\gamma N^{-1/2},\theta}-\mu_N\rrVert  _{\mathrm{TV}}$ converges
uniformly in $\gamma\in[0,K]$ and in $\theta$, to
%
\begin{equation}
\frac{1}{\sqrt{8\pi}}\int\bigl\llvert e^{\gamma u-\gamma^2/2}-1\bigr\rrvert
e^{-\sfrac
{u^2}{2}}\,\DD u.
\end{equation}
The conclusion then follows by performing a few changes of variables.

\begin{pf*}{Proof of Lemma~\ref{GTLC}}
Statement (ii) follows from a more general statement on concentration
for Lipshitz functional on $\gO_N$: Proposition~\ref{lipitz} is
proved in the \hyperref[techos]{Appendix}.

To ensure that the convergence holds uniformly in $\theta$, the reader
can check that
all the bounds present in the proof do not depend on $\theta$.
In the remainder of the paper, we will use the notation
%
\begin{equation}
\label{bsin} \bsin(x)=\bsin_{\theta}(x):= \sin \biggl(
\frac{x\pi}{N}+\theta \biggr).
\end{equation}
As most computations do not depend on the value of $\theta$, we omit
it from the notation most of the time in the remainder of the paper.
Note that $a(\eta)$ satisfies trivially
$\mu_N(a(\eta))=0$.
Let us show that the
variance is asymptotically equivalent to $N$.
%
\begin{equation}
\mu_N\bigl(a(\eta)^2\bigr)=\sum
_{x\in\bbZ_{2N}} \bsin(x)^2+\mathop{\sum
_{(x,y)\in\bbZ_{2N}}}_{x\ne y}\bsin(x)\bsin(y) \bbE\bigl[\eta(x)
\eta(y)\bigr].
\end{equation}
The first term is equal to $N$. As for the second term, as we have
$\bbE[\eta(x)\eta(y)]= -1/(2N+1)$ it is equal to
%
\begin{equation}
\frac{1} {(2N-1)}\sum_{x\in\bbZ_{2N}}\bsin(x)^2=
\frac{N}{2N-1}.
\end{equation}
To show the convergence to a Gaussian variable, we will use the
martingale central limit theorem \cite{cfBrown}, Theorem 1.
Let $(M^N_i)_{i=0}^{2N-1}$ be the martingale defined by
%
\begin{equation}
M^N_i:= \mu_N \bigl( a(\eta) \mid \bigl(
\eta(x)\bigr)_{x=1}^i \bigr).
\end{equation}
It satisfies $M^N_0=0$ and $M^N_{2N-1}=a(\eta)$.
Set
%
\begin{equation}
\Delta M_i:= M^{N}_{i+1}-M^{N}_i
\end{equation}
and
%
\begin{equation}
\sigma^2_{i,N}= \mu_N \bigl( (\Delta
M_i)^2 \mid \bigl(\eta(x)\bigr)_{x=1}^i
\bigr).
\end{equation}
To apply the central limit theorem the martingale $M_i$,
one must only check that
%
\begin{equation}
\label{sigma2N} \sigma^2_N:=\frac{1}{N}\sum
_{i=0}^{2N-2}\sigma^2_{i,N},
\end{equation}
converges to one, in probability (there are in fact other assumptions
to check; see~\cite{cfBrown} but in our case they are trivially satisfied).

For $A\subset\bbZ_{2N}$, we let $\eta(A)$ denote the number of
particles in the set $A$,
%
\begin{equation}
\label{distrib} \eta(A):=\sum_{x\in A}
\ind_{\{\eta_{x}=1\}}.
\end{equation}
Let us fix $i$ and set $k=k(\eta,i):=\eta([1,i])$.
A simple computation gives
%
\begin{equation}
\Delta M_i= \cases{
\displaystyle\frac{N+k-i}{2N-i} \Biggl(\bsin(i+1)-
\frac{1}{2N-i-1}\sum_{x=i+2}^{2N}\bsin(x)
\Biggr),
\cr
\qquad\mbox{if $\eta(i+1)=1$,}
\vspace*{3pt}\cr
\displaystyle\frac{N-k}{2N-i} \Biggl(\bsin(i+1)-
\frac{1}{2N-i-1}\sum_{x=i+2}^{2N}\bsin(x)
\Biggr),
\cr
\qquad\mbox{if $\eta(i+1)=-1$.}}
\end{equation}
As the first and second option in \eqref{distrib} have respective
probability $(N-k)/(2N-i)$ and $(N+k-i)/(2N-i)$,
we have
%
\begin{equation}
\sigma^2_{i,N}= \frac{2(N+k-i)(N-k)}{(2N-i)^2} \Biggl( \bsin(i+1)-
\frac{1}{2N-i-1}\sum_{x=i+2}^{2N} \bsin(x)
\Biggr)^2.\hspace*{-30pt}
\end{equation}
Once this is done, we just need to check the following facts to conclude:
\begin{longlist}[(a)]
\item[(a)] For all $i$, $\sigma^2_{i,N}$ is almost surely smaller
than $8$.
\item[(b)] For all $N$ sufficiently large, for all $i\in[0,2N-\sqrt
{N}]$ we have
%
\begin{eqnarray}
&& \mu_N \Biggl( \Biggl\llvert \sigma^2_{i,N}-
\frac{1}{2} \Biggl( \bsin(i+1)- \frac{1}{2N-i-1}\sum
_{x=i+2}^{2N} \bsin(x) \Biggr)^2\Biggr
\rrvert \ge N^{-1/20} \Biggr)
\nonumber\\[-8pt]\\[-8pt]\nonumber
&&\qquad \le N^{-1/20}.
\end{eqnarray}
\item[(c)] We have the following convergence:
%
\begin{equation}
\lim_{N\to\infty} \frac{1}{N}\sum
_{i=1}^{2N-1} \Biggl( \bsin(i+1)- \frac{1}{2N-i-1}
\sum_{x=i+2}^{2N} \bsin(x)
\Biggr)^2=2.
\end{equation}
\end{longlist}
From the these three claims, it is rather standard to show that
$\sigma^2_{N}$ converges to $1$ in probability and we leave it as an
exercise to the reader.
Item~(a) is obvious, item (b) follows from computing the mean and
variance of $k(\eta,i)$ [which are, resp., equal to
$i/2$ and $i(2N-i)/4(2N-1)$] and applying the Markov inequality.
As for $(c)$, it can be checked via a tedious but simple computation.
\end{pf*}

\section{\texorpdfstring{Decomposing the proof of Proposition \protect\ref{mainresult}.}{Decomposing the proof of Proposition 2.3}}\label{decompo}

We present in this section the main steps of the proof of our main
result and the heuristics behind it.

\subsection{\texorpdfstring{Why coupling with $\nu^{N,\alpha,\theta}$?}{Why coupling with nu{N,alpha,theta}?}}

First, let us try to understand why $\nu^{N,\alpha,\theta}$ gives a
good approximation of the $P^\chi_t$, via an analysis of the particle
density and fluctuation.
Let
%
\begin{equation}
u^\chi(x,t):= \bbE\bigl[\eta^\chi_t(x)\bigr]
\end{equation}
denote the expected particle density (with our convention it can be
negative since empty sites count for $-1$).
It is standard to check that $u^\chi$ is the solution of the discrete
heat-equation
%
\begin{equation}
\label{disheat} %
\cases{
\partial_t u(x,t):= \gD u(x,t),
\cr
u(x,0):=\chi(x),}
\end{equation}
where $\gD$ denotes the discrete Laplacian
%
\begin{equation}
\gD u(x,t):= u(x+1,t)+u(x-1,t)-2u(x,t).
\end{equation}
Here and in what follows if $f$ is a function of $\bbZ_{2N}$
(identified to a periodic function of $\bbZ$) such that
%
\begin{equation}
\sum_{x\in\bbZ_{2N}} f(x)=0,
\end{equation}
and $x$ and $y$ are two elements of $\bbZ_{2N}$ and
$x_0\le y_0$ two elements of $\bbZ$ which are, respectively, equal to
$x$ and $y$ modulo $2N$.
Then we use the notation
$
\sum_{z=x}^y f(z)$,
to denote the sum
$
\sum_{z=x_0}^{y_0} f(z)$.
It is straightforward that it does not depend on the particular choice
of $x_0$ and $y_0$ once $x$ and $y$ are fixed.
Let us quickly investigate the fluctuations of the integrated density
of particle
%
\begin{equation}
\xi(\eta) (x):=\sum_{z=1}^x \eta(z).
\end{equation}
At equilibrium, $\xi(\eta)$ is a simple random-walk conditioned to
return to zero after $2N$ steps.
It is centered and has Gaussian fluctuations of order $\sqrt{N}$.
In \cite{cfLac2}, we have proved that
the fluctuation of $\xi(\eta^\chi_t)(x)$ around its mean [given by
$\sum_{z=1}^x u^{\chi}(z,t)$] are in fact always of order $\sqrt{N}$.

This gives the intuition that much before the mixing time, the law of
$\eta^\chi_t$ might approximately be that
of $2N$ independent $\pm1$ Bernoulli variables,\vspace*{1pt} each with bias
$u^{\chi}(x,t)$, conditioned to $\sum_{x\in\bbZ_{2N}} \eta^\chi_t(x)=0$.

For\vspace*{2pt} $t\ge\frac{N^2}{4\pi^2}\log N$, $u^{\chi}(x)$ is very
well approximated by a sinusoid function (see Lemma~\ref{sinso})
%
\begin{equation}
u^{\chi}(x, t)\approx b(\chi)e^{-\gl_N t}\sin \biggl(
\frac{\pi
x}{N}+\theta(\chi) \biggr),
\end{equation}
and the conditioned law of independent Bernoulli described above is
very close in total variation to
$\nu^{N,b(\chi)e^{-\gl_N t},\theta(\chi)}$.

\subsection{\texorpdfstring{How to do it.}{How to do it}}

Let us first write here the rigorous result concerning the fluctuation
around the expected density of particle.

\begin{proposition}\label{sinusoid}
There exists a constant $c>0$ such that for all $N$ sufficiently large,
for all $t\ge\frac{1}{4}(\gl_N)^{-1}\log N$, we have
%
\begin{eqnarray}
&& P^{\chi}_t \Biggl[ \exists x, y \in\bbZ_{2N},
\Biggl\llvert \sum_{z=x+1}^y \biggl[ \eta(z)
- e^{-\gl_{N}t} b(\chi) \sin \biggl(\frac{\pi z}{N}+\theta(\chi) \biggr)
\biggr] \Biggr\rrvert \ge s \sqrt{N} \Biggr]
\nonumber\\[-8pt]\\[-8pt]\nonumber
&&\qquad \le2e^{-cs^2}
\end{eqnarray}
\end{proposition}

In particular, we know that with high probability, $\eta_t^\chi$ lies
in the set
%
\begin{eqnarray}
\mathcal G^N_{\alpha,\theta} &:=& \Biggl\{ \eta\in\gO_N
\Big| \max_{x, y \in\bbZ_{2N}} \Biggl\llvert \sum
_{z=x+1}^y \biggl(\eta(z)
\nonumber\\[-8pt]\\[-8pt]\nonumber
&&{} - \alpha\sin \biggl(
\frac{\pi
z}{N}+\theta \biggr) \biggr) \Biggr\rrvert \le\sqrt{N} \log\log N
\Biggr\}
\end{eqnarray}
with $\alpha$ and $\theta$ being chosen, respectively, equal to
$e^{-\gl_{N}t} b(\chi)$ and $\theta(\chi)$.

To prove Proposition~\ref{mainresult}, it is sufficient to prove that:
\begin{itemize}
\item within a time
$N^2(\log N)^{1/2}$ (i.e., a shorter time-scale than the mixing time),
one can couple
a dynamics with initial condition $\chi\in\mathcal G^N_{\alpha,\theta}$, with one
with initial condition distributed like $\nu^{N,\alpha,\theta}$.
\item
the family of measure $(\nu^{N,\alpha, \theta})$ is almost preserved
by the dynamics
in the sense that applying the semi-group $P_t$ to it only changes the
value of $\alpha$.
\end{itemize}
Both of these statements hold provided $\alpha$ is sufficiently small,
and are stated as two propositions below.
More\vspace*{1pt} precisely
Let $\nu^{N,\alpha, \theta}_t$ be the law of a system started with
initial configuration
$\nu^{N,\alpha,\theta}$
%
\begin{equation}
\nu^{N,\alpha, \theta}_t(\eta):=\sum_{\eta'\in\gO_N}
\nu ^{N,\alpha, \theta}_t P_t\bigl(\eta',\eta
\bigr).
\end{equation}

\begin{proposition}\label{prop1}
For all $N$ sufficiently large,
for all $\theta\in[0,2\pi)$, for all $\alpha\le2N^{-3/7}$, we have
for all $\chi\in\mathcal G_{\alpha,\theta}$, for all $t\ge N^2
(\log N)^{1/2}$
%
\begin{equation}
\bigl\llVert P^\chi_t- \nu^{N,\alpha, \theta}_t
\bigr\rrVert \le\frac{1}{2\log\log N}.
\end{equation}
\end{proposition}

\begin{proposition}\label{densitevol}
There exists a constant $C$ such that for all $N$ and all \mbox{$\alpha>0$},
%
\begin{equation}
\sup_{t\ge0}\bigl\llVert \nu^{N,\alpha,\theta}_{t}-
\nu^{N, \alpha e^{-\gl
_{N}t},\theta}\bigr\rrVert _{\mathrm{TV}}\le C\alpha^2 \sqrt{N}.
\end{equation}
\end{proposition}

\begin{pf*}{Proof of Proposition~\ref{mainresult}}
We have for any $\chi$ in $\beta$ for
$t\ge t_0:= 3/7(\gl_N)^{-1}$ we have, using the triangular inequality
%
\begin{eqnarray}
&& \bigl\llVert P^\chi_t- \nu^{N,b(\chi)e^{-\gl_N t}, \theta(\chi)} \bigr\rrVert
_{\mathrm{TV}} \nonumber
\\
&&\qquad \le\sum_{\chi' \in\gO_N} P_t\bigl(
\chi,\chi'\bigr) \bigl\llVert P^{\chi'}_{t-t_0}-
\nu^{N,b(\chi)e^{-\gl_N t}, \theta(\chi)}\bigr\rrVert _{\mathrm{TV}}
\nonumber\\[-8pt]\\[-8pt]\nonumber
&&\qquad \le P^\chi_{t_0} \bigl(\eta\notin\cG^N_{b(\chi)N^{-3/7},\theta
(\chi)} \bigr)
\\
&&\quad\qquad{} +\max_{\chi'\in\cG^N_{b(\chi)N^{-3/7},\theta
(\chi)}} \bigl\llVert P^{\chi'}_{t-t_0}-
\nu^{N,b(\chi)e^{-\gl_Nt}, \theta(\chi)}\bigr\rrVert _{\mathrm{TV}}.\nonumber
\end{eqnarray}
According to Proposition~\ref{sinusoid}, we have
%
\begin{equation}
P^\chi_{t_0} \bigl(\eta\notin\cG^N_{b(\chi)N^{-3/7},\theta(\chi
)}
\bigr)\le\frac{1}{\log N}.
\end{equation}
We note that for $\chi'\in\cG^N_{b(\chi)N^{-3/7},\theta(\chi)}$
we have
%
\begin{eqnarray}
&& \bigl\llVert P^{\chi'}_{t-t_0}-\nu^{N,b(\chi)e^{-\gl_N t}, \theta(\chi)}\bigr\rrVert
_{\mathrm{TV}}\nonumber
\\
&&\qquad \le \bigl\llVert P^{\chi'}_{t-t_0}-
\nu^{N,b(\chi)N^{-3/7}, \theta(\chi)}_{t-t_0}\bigr\rrVert _{\mathrm{TV}}
\\
&&\quad\qquad{} + \bigl\llVert \nu^{N,b(\chi)N^{-3/7}, \theta(\chi)}_{t-t_0}-\nu^{N,b(\chi
)e^{-\gl_N t}, \theta(\chi)}\bigr\rrVert
_{\mathrm{TV}}.\nonumber
\end{eqnarray}
The first term is smaller than $\frac{1}{2\log\log N }$ according to
Proposition~\ref{prop1} as soon as
\[
t\ge t_0+N^2\sqrt{\log N}.
\]
Proposition~\ref{densitevol} ensures that the second term is smaller than
$(\log N)^{-1}$, hence we can conclude.
\end{pf*}

\section{\texorpdfstring{Proof of Proposition \protect\ref{sinusoid}.}{Proof of Proposition 4.1}}\label{fluctutu}

This statement is in fact mostly derived from the statement about
fluctuations proved in \cite{cfLac2}
which we state now.

\begin{proposition}[(\cite{cfLac2}, Proposition 4.1)]\label{sinusoid2}
There exists a constant $c>0$ such that for all $t\ge0$, for all $s\ge
0$, for all $\chi\in\gO_N$ we have
%
\begin{equation}
P^\chi_t \Biggl[ \exists x, y \in\bbZ_{2N},
\Biggl\llvert \sum_{z=x+1}^y \bigl(\eta(z)-
u^\chi(z,t) \bigr) \Biggr\rrvert \ge s \sqrt{N} \Biggr]
\le2e^{-cs^2}.
\end{equation}
\end{proposition}

\begin{rem}
Note that in \cite{cfLac2}, $t\ge3N^2$ is required (that would be in
fact $t\ge12 N^2$ in our setup because we work on $\bbZ_{2N}$ instead
of $\bbZ_N$), but this is only to treat the case of an arbitrary
number of particles. The reader can check from the proof that this
assumption is only needed to check
\cite{cfLac2}, equation~(4.4), which is obviously valid for all $t\ge
0$ when we have $N$ particles on $2N$ sites.
\end{rem}

With this result, Proposition~\ref{sinusoid} follows from a basic
analysis of the Fourier coefficients of the solution of \eqref{disheat}.

\begin{lemma}\label{sinso}
For all $t\ge\frac{1}{4}(\gl_N)^{-1}\log N$, we have
%
\begin{equation}
\max_{x\in\bbZ_{2N}} \biggl\llvert u^{\chi}(x,t)-
e^{-\gl_{N}t} b(\chi) \sin \biggl(\frac{\pi x}{N}+\theta(\chi) \biggr)
\biggr\rrvert \le4 N^{-1/2}.
\end{equation}
\end{lemma}

\begin{pf}
Let us decompose $u^\chi$ along its Fourier modes for the
heat-equation.

As in Section~\ref{intro}, we prefer to have only one coefficient per
eigenspace, and thus,
for $i=2,\dots,N-1$, introduce $\theta_i(\chi)$ to be
the phase of the projection.
We let $\theta_i(\chi)$ be either the unique solution of
%
\begin{eqnarray}
\label{deftheta} \sum_{x\in\bbZ_{2N}} \chi(x) \cos \biggl(
\frac{i\pi x}{N}+\theta \biggr) &=&0,
\nonumber\\[-8pt]\\[-8pt]\nonumber
\sum_{x\in\bbZ_{2N}} \chi(x) \sin \biggl(\frac{i\pi x}{N}+
\theta \biggr) &>&0
\end{eqnarray}
or $0$ if all $\theta$ are solution, and let us set
%
\begin{eqnarray}
\label{defb} b_i(\chi)&:=&\frac{1}{N}\sum
_{z\in\bbZ_{2N} }\chi(x) \sin \biggl(\frac{i\pi x}{N}+
\theta_i(\chi) \biggr),
\nonumber\\[-8pt]\\[-8pt]\nonumber
b_{N}(\chi)&:=&\frac{1}{2N} \sum_{z\in\bbZ_{2N} }(-1)^{\llvert  x\rrvert  }
\chi(x).
\end{eqnarray}
We have
%
\begin{equation}
\chi(x):=\sum_{i=1}^{N-1} b_i(
\chi) \sin \biggl(\frac{i\pi
x}{N}+\theta_i(\chi)
\biggr)+b_N(\chi) (-1)^{\llvert  x\rrvert  }.
\end{equation}
As the functions
$x \mapsto\sin (\frac{i\pi x}{N}+\theta_i )$ are
eigenfunctions of the discrete Laplacian with respective eigenvalues,
%
\begin{equation}
-\gl_{i,N}:=2 \biggl(1-\cos \biggl(\frac{i\pi}{N} \biggr) \biggr),
\end{equation}
we have for all $t\ge0$,
%
\begin{equation}
u(x,t)= e^{-\gl_{N}t} b(\chi) \sin \biggl(\frac{\pi x}{N}+\theta (\chi)
\biggr)+R(\chi,t,x),
\end{equation}
where
%
\begin{equation}
R(\chi,t,x):= \sum_{i=2}^{N-1}
e^{-\gl_{i,N}t} b_i(\chi) \sin \biggl(\frac{i\pi
x}{N}+
\theta_i(\chi) \biggr)+ b_N(\chi)e^{-2t}.
\end{equation}
Noticing that $\llvert   b_i(\chi)\rrvert   \le2$ and that for all $N\ge2$,
%
\begin{equation}
\forall i\in\{2,\dots,N\}, \qquad\gl_{i,N}\ge i \gl_N,
\end{equation}
we have for all $t\ge(\gl_N)^{-1}$
%
\begin{equation}
\bigl\llvert R(\chi,t,x)\bigr\rrvert \le2 \sum_{i=2}^N
e^{-i \gl_N t}= \frac{ 2e^{-2 \gl
_N t}}{1-e^{-\gl_N t}} \le4 e^{-2\gl_N t}.
\end{equation}
Hence, we have the result.
\end{pf}

\section{\texorpdfstring{Proof of Proposition \protect\ref{prop1}.}{Proof of Proposition 4.2}}\label{coupel}
We assume without loss of generality that $\alpha$ is nonnegative, and write $\nu^\alpha$
and $\nu^{\alpha}_t$ for $\nu^{N,\alpha,0}$ and
$\nu^{N,\alpha,0}_t$, and $a(\eta)$ for $a_0(\eta)$.

\subsection{\texorpdfstring{Properties of $\nu^\alpha$.}{Properties of nualpha}}\label{repair}

In this section, we check several properties for $\nu^{\alpha}$.
While the results are quite intuitive,
their proof is quite technical and we have decided to postpone them to
Appendix~\ref{techos}.
First, we want to ensure that it has the right density of the particle.

\begin{proposition}\label{rightdensity}
There exists a constant $C$ such that for all $\alpha\le1$ we have
%
\begin{equation}
\label{dsss} \sup_{x\in\bbZ_{2N}}\bigl\llvert \nu^{\alpha}\bigl(
\eta(x)\bigr)-\alpha\bsin(x)\bigr\rrvert \le C\bigl(\alpha^2+N^{-2}
\bigr).
\end{equation}
\end{proposition}

Then we have to check that the fluctuations are not larger than
$\sqrt{N}$.

\begin{proposition}\label{fluctouse}
There exists constant $c$ such that for all $N>0$, for all $\llvert   \alpha \rrvert  \le N^{-1/4}$, and $t\ge0$
%
\begin{equation}
\label{fructoe} \nu^\alpha_t \Biggl[ \exists x, y \in
\bbZ_{2N}, \Biggl\llvert \sum_{z=x+1}^y
\bigl(\eta(z)- \alpha e^{-\gl_N t}\bsin(z) \bigr) \Biggr\rrvert \ge s
\sqrt{N} \Biggr]\le2e^{-cs^2}.
\end{equation}
\end{proposition}

Finally, we want to check that if one starts from distribution $\nu
^{\alpha}$
there is a positive density of sites where $\eta(z)\ne\eta(z+1)$,
that is, of locations where jumps of the particle can occur.
The utility of such a statement will be become clear in the next
section when we construct the dynamical coupling.
For a probability measure $\nu$ defined on $\gO_N$, we let $\bbP
^{\nu}$ be the law of the Markov chain $(\eta_t)_{t\ge0}$
starting from $\eta_0$ distributed like $\nu$.
Set
%
\begin{equation}
j(x,y,\eta):= \bigl\{ z\in[x,y] \mid \eta(z)\ne\eta(z+1) \bigr\}
\end{equation}
and
%
\begin{equation}
\label{defce} \qquad\cE:= \bigl\{\eta\in\gO_N \mid \forall(x,y) \in
\bbZ^2_{2N}, \# [x,y]\ge N^{1/4}\Rightarrow j(x,y,
\eta) \ge\tfrac{1}{4}\#[x,y] \bigr\}.
\end{equation}

\begin{proposition}\label{iamzea}
There exist a constant $c$ such that
for $N$ sufficiently large, for all $\llvert   \alpha \rrvert  \le c N^{-3/8} $
%
\begin{equation}
\bbP^{\nu^{\alpha}}\bigl[ \exists t\le N^3, \eta_t
\notin\cE\bigr]\le e^{-cN^{1/4}}.
\end{equation}
\end{proposition}

\begin{rem}
The power exponents for $N$ in Proposition \ref{iamzea} are rather arbitrary and other choices would also fit.
The important result is that the probability tends to zero.
\end{rem}

\subsection{\texorpdfstring{The $\xi$ dynamics.}{The xi dynamics}}

We introduce in this section an auxiliary dynamics (the same as in
\cite{cfLac2}) which is used to couple
$P_t^\chi$ with $\chi\in\mathcal G_{\alpha}$ (we use this notation
for $\mathcal G^N_{\alpha,0}$) with $\nu^\alpha_t$.
The idea of using interface dynamics to study particle system dates is
not new and is already present in
the seminal paper of Rost about the asymmetric exclusion on the line
\cite{cfRost}
(for the use of this technique for mixing time related issues, see
\cite{cfWilson,cfLac1,cfLac2}).
In \cite{cfWilson,cfLac1}, the height function is introduced mainly
to have a better intuition on an order which can be defined without
the interface representation. Let us stress that here, on the contrary,
the interface dynamics is used to perform a monotone coupling that
could not be constructed by considering only the original chain.

Let us consider the set of discrete height functions of the circle.
%
\begin{equation}
\gO'_{N}:= \bigl\{\xi: \bbZ_{2N} \to\bbZ \mid
\xi(0)\in2\bbZ, \forall x\in\bbZ_{2N}, \bigl\llvert \xi(x)-\xi(x+1)
\bigr\rrvert =1 \bigr\}.
\end{equation}
Given $\xi$ in $\gO'_{N}$, we define $\xi^x$ as
%
\begin{equation}
\cases{ \xi^x(y)=\xi(y), &\quad$\forall y \ne x$,
\cr
\xi^x(x)=\xi(x+1)+\xi(x-1)-2\xi(x).}
\end{equation}
We let $\xi_t$ be the irreducible Markov chain on $\gO'_{N}$ whose
transition rates $p$ are given by
%
\begin{equation}
\label{cromik} \cases{ p\bigl(\xi,\xi^x\bigr)=1, &\quad$\forall x\in
\bbZ_N$,
\cr
p\bigl(\xi,\xi'\bigr)=0, &\quad if $
\xi'\notin\bigl\{ \xi^x \mid x\in\bbZ_N\bigr
\}$.}
\end{equation}
We call this dynamics the corner-flip dynamics, as the transition $\xi
\to\xi^x$ corresponds to flipping
either a local maximum of $\xi$ (a ``corner'' for the graph of $\xi$)
to a local minimum {e vice versa}.
It is of course not positive recurrent, as the state space is infinite
and translation invariant for the dynamics, however,
it is irreducible and recurrent.

\begin{figure}

\includegraphics{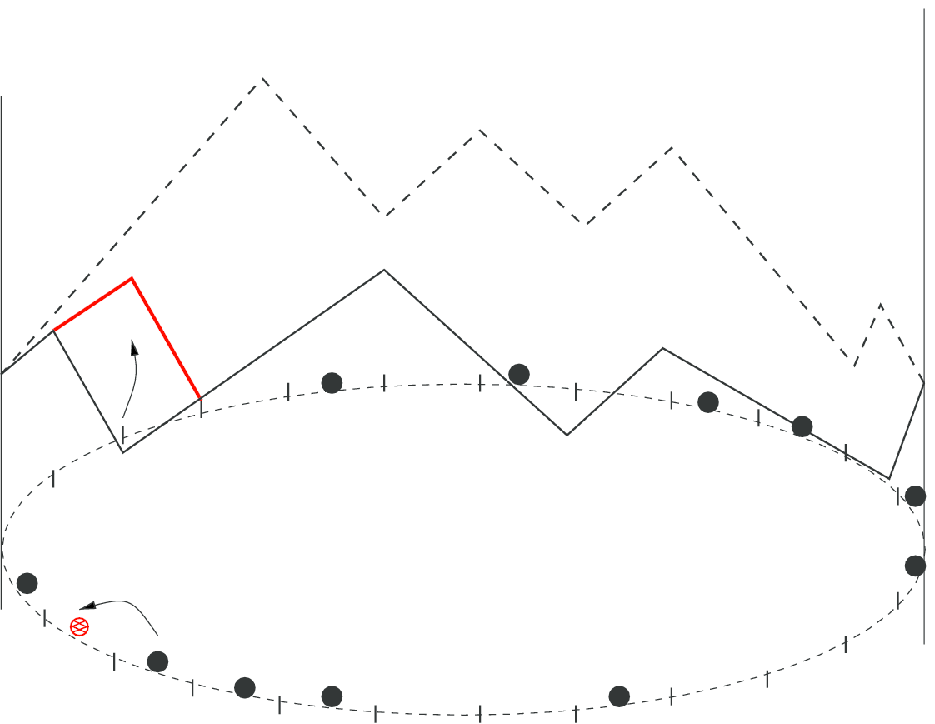}

\caption{The correspondence between the exclusion process and the
corner-flip dynamics.
A particle jump and its corner-flip counterpart are indicated by arrows.
Note that this is not a one-to-one mapping
as a particle configuration gives the height function only modulo translation.}
\label{partisys}
\end{figure}

The reader can check that $\gO'_{N}$ is mapped onto $\gO_{N}$, by the
transformation
$\xi\mapsto\grad\xi$ where
%
\begin{equation}
\label{defgrad} \nabla\xi(x):=\xi(x+1)-\xi(x)
\end{equation}
and that
the image of the corner-flip dynamics $(\nabla\xi_t)_{t\ge0}$ is the
simple exclusion process (see Figure~\ref{partisys}).

There is a natural order on the set $\gO'_{N}$
defined by
%
\begin{equation}
\label{deforder} \xi\ge\xi' \quad\Leftrightarrow\quad\forall x\in
\bbZ_{2N},\qquad \xi (x)\ge\xi'(x),
\end{equation}
and we can construct a \textit{grand coupling} for the Markov chain
which preserves this order.

\subsection{\texorpdfstring{The graphical construction.}{The graphical construction}}\label{grafff}

We introduce in this section an order preserving grand-coupling on $\gO'_{N}$.
For $\zeta\in\gO'_N$, $(\xi^{\zeta}_t)_{t\ge0}$ denotes the
Markov chain with initial condition $\zeta$.
We want to construct all the $(\xi^{\zeta}_t)_{t\ge0}$ on a same
probability space in a way that
%
\begin{equation}
\label{preserving} \forall\zeta, \zeta' \in\gO'_N,
\bigl(\zeta\ge\zeta'\bigr) \Rightarrow \bigl( \forall t\ge0,
\xi^{\zeta}_t\ge\xi^{\zeta'}_t \bigr).
\end{equation}
Of course, there are several options for such a grand coupling. We want
to choose one which is such
that, eventually, the trajectories starting from different initial
conditions coalesce almost surely (at a random time)
%
\begin{equation}
\forall\zeta, \zeta' \in\gO'_N, \exists
T_{\zeta,\zeta
'}<\infty, \forall t\ge T_{\zeta,\zeta'}, \qquad
\xi^{\zeta
}_t=\xi^{\zeta'}_t.
\end{equation}
Of course, we want the coalescing time to be as short as possible.
To reach this aim, we make the different corner flips for different
trajectories as independent as can be while still satisfying
\eqref{preserving}.

Let us present the construction.
The evolution of the $(\xi_t)_{t\ge0}$ is completely determined by
auxiliary Poisson processes which we call clock processes.
Set
\[
\Theta:= \bigl\{ (x,z) \mid x\in\bbZ_N\mbox{ and } z\in 2
\bbZ+(-1)^x \bigr\}.
\]

And set $\cT^\uparrow$ and $\cT^\downarrow$ to be two independent
rate-one clock processes
indexed by $\Theta$ ($\cT^\uparrow_{\omega}$ and
$\cT^\downarrow_{\omega}$ are two independent Poisson processes of
intensity one of each $\omega\in\Theta$).
The trajectory of $\xi_t$ given $(\cT^\uparrow,\cT^\downarrow)$ is
given by the following construction:
\begin{itemize}
\item$\xi_t$ is a c\`adl\`ag, and does not jump until one of the
clocks indexed by $(x,\xi_t(x))$, $x\in\bbZ_{2N}$ rings.
\item If $\cT^\downarrow_{(x,\xi_{t^-}(x))}$ rings at time $t$ and
$x$ is a local maximum for $\xi_{t^-}$, then
$\xi_{t}=\xi^x_{t^-}$.
\item If $\cT^\uparrow_{(x,\xi_{t^-}(x))}$ rings at time $t$ and
$x$ is a local minimum for $\xi_{t^-}$, then
$\xi_{t}=\xi^x_{t^-}$.
\end{itemize}

\subsection{\texorpdfstring{Construction the initial condition for $\xi^0$, $\xi^1$ and $\xi^2$.}{Construction the initial condition for xi0, xi1 and xi2}}

Given $\chi\in\cG_\alpha$, we let $(\xi^0_t)$ the trajectory of
the Markov chain with transitions rates
\eqref{cromik} starting from initial condition
%
\begin{equation}
\label{fryied} \xi_{0}^0(x):=\sum
_{z=0}^x \chi(x).
\end{equation}
Note that for all $t\ge0$ we have
%
\begin{equation}
\label{laloi} \bbP \bigl[ \grad\xi^0_t\in\cdot
\bigr]=P^\chi_{t}.
\end{equation}

Our idea is to construct another dynamic $\xi^1_t$ which starts with
$\grad\xi^1_0$ distributed like $\nu^{\alpha}$
which coalesces with $\xi^0_t$ within time $N^2\sqrt{\log N}$.
In fact, it turns out more practical to define not one but two dynamics
$\xi^1$ and $\xi^2$ to couple with $\xi^0$.
We let $\bbP$ denote the law of $(\xi^0_t,\xi^1_t,\xi^2_t)_{t\ge
0}$, and we impose
%
\begin{equation}
\label{equistart0} \bbP \bigl[\grad\xi^1_0\in\cdot \bigr]=
\bbP \bigl[\grad\xi ^2_0\in\cdot \bigr]=
\nu^\alpha.
\end{equation}
Note that this implies for all $t\ge0$
%
\begin{equation}
\label{equistart} \bbP \bigl[\grad\xi^1_t\in\cdot \bigr]=
\bbP \bigl[\grad\xi ^2_t\in\cdot \bigr]=
\nu^\alpha_t.
\end{equation}
We impose also the condition
%
\begin{equation}
\label{rorderring0} \xi^1_0\le\xi^0_0
\le\xi^2_0,
\end{equation}
and use the graphical coupling introduced in the previous section to
construct the trajectory
of $(\xi^i_t)_{t\ge0}$, $i=0,1,2$. Hence, the order is conserved at
all time
%
\begin{equation}
\label{rorderring} \forall t\ge0,\qquad\xi^1_t\le
\xi^0_t\le\xi^2_t.
\end{equation}
Let us now explain our construction of the initial conditions.
We start with $\eta_0$ distributed like $\nu^{\alpha}$ and we will
choose $\xi^1_0$ and $\xi^2_0$ such that
%
\begin{equation}
\grad\xi^1_0=\grad\xi^2_0=
\eta_0.
\end{equation}
We set for arbitrary $\eta\in\gO_{N}$, or $\xi\in\gO'_{N}$
%
\begin{eqnarray}
\label{defh} H_{t,\alpha}(\eta)&:=&\max_{x,y\in\bbZ_N} \Biggl
\llvert \sum_{z=x+1}^y \eta(z)-e^{-\gl_Nt}
\bsin(t)\Biggr\rrvert,
\nonumber\\[-8pt]\\[-8pt]\nonumber
H_{t,\alpha}(\xi)&:=& H_{t,\alpha}(\grad\xi).
\end{eqnarray}
We also set
%
\begin{equation}
\label{defh0} \mathcal H_0:=2 \bigl\lceil\bigl(H_{0,\alpha}(
\eta_0)+\sqrt{N}\log\log N\bigr) /2 \bigr\rceil
\end{equation}
and
%
\begin{eqnarray}
\xi^1_0(x)&:=&\sum_{z=1}^x
\eta_0(z)-\cH_0,
\nonumber\\[-8pt]\\[-8pt]\nonumber
\xi^2_0(x)&:=&\sum_{z=1}^x
\eta_0(z)+\cH_0.
\end{eqnarray}
The fact that \eqref{rorderring0} is satisfied follows from the
definition of $\mathcal G_\alpha$ and
that of $\cH_0$.
Note also that from Proposition~\ref{fluctouse} applied at $t=0$, we have
%
\begin{equation}
\bbP [\mathcal H_0\ge2\sqrt{N}\log\log N ]\le(\log
N)^{-1}.
\end{equation}
To prove Proposition~\ref{prop1}, it is sufficient to prove that $\xi
^1_t$ and $\xi^2_t$
typically coalesce within a time $N^2\sqrt{\log N}$. More precisely,
we have the following.

\begin{proposition}\label{nahnou}
For sufficiently large $N$, for all $\alpha\le2N^{-3/7} $
for $(\xi^1_t)_{t\le0}$, $(\xi^2_t)_{t\ge0}$, constructed as above,
we have
%
\begin{equation}
\label{nahnou1} \bbP \bigl[ \xi^1_{N^2\sqrt{\log N} }\ne
\xi^2_{N^2\sqrt{\log N}} \bigr]\le\frac{1}{2\log\log N}.
\end{equation}
\end{proposition}

Proposition~\ref{nahnou} is proved in Sections~\ref{croook} and~\ref{craak}.

\begin{pf*}{Proof of Proposition~\ref{prop1}}
Let $\chi$ in $\mathcal G_\alpha$ be fixed and consider the dynamics
$\xi^i$, $i=0,1,2$ constructed as above.
From \eqref{rorderring}, we have
%
\begin{equation}
\xi^1_t=\xi^2_t \quad\Rightarrow
\quad\xi^1_t=\xi^0_t.
\end{equation}
Recalling \eqref{laloi} and \eqref{equistart}, we have for any $t>0$
%
\begin{equation}
\bigl\llVert P^\chi_{t}-\nu^\alpha_t
\bigr\rrVert _{\mathrm{TV}}\le\bbP \bigl[ \grad\xi^0_{t}
\ne \grad\xi^1_{t} \bigr]\le\bbP \bigl[
\xi^0_{t}\ne\xi ^1_{t} \bigr] \le
\bbP \bigl[\xi^1_{t}\ne\xi^2_{t}
\bigr].
\end{equation}
Hence, Proposition~\ref{nahnou} implies the result.
\end{pf*}

\subsection{\texorpdfstring{The randomly walking area.}{The randomly walking area}}\label{croook}

Let us set
%
\begin{equation}
A(t)=\frac{1}{2}\sum_{x\in\bbZ_2}
\xi^2_t(x)-\xi^1_t(x).
\end{equation}
The reader can check that $A(t)$ is an integer.
Because of \eqref{rorderring}, we remark that $A(t)$ is always
positive, and hence that
$\xi^1$ and $\xi^2$ merge at time
%
\begin{equation}
\label{deftau} \tau:=\inf\bigl\{ t \ge0 \mid A(t)=0 \bigr\}.
\end{equation}
As $A(t)$
is an integer valued martingale which only makes $\pm1$ jumps, it is
to be a time changed symmetric nearest neighbor walk on $\bbZ_+$.
In order to get a bound for
\[
\bbP[\tau\le t],
\]
we need to have a reasonable control over the time change, that is, the
jump rate of $A(t)$.
It depends on
the particular configuration $(\xi^1_t,\xi^2_t)$ the system sits on:
it is given by
the number of places where corners can flip independently for $\xi
^1_t$ and $\xi^2_t$.
More precisely,
set
%
\begin{eqnarray}
\label{defui} U_i(t)&:=&\bigl\{x\in\bbZ_N \mid
\xi^i_t\mbox{ has a local extremum at } x \mbox{ and}
\nonumber\\[-8pt]\\[-8pt]\nonumber
&&{} \exists y\in\{x-1,x,x+1\}, \xi^2_t(y)>
\xi^1_t(y)\bigr\}.
\end{eqnarray}
The jump rate of $A(t)$ is given by
%
\begin{equation}
u(t):=\#U_1(t)+\#U_2(t).
\end{equation}
For $t\le\int_0^\tau u(t) \,\DD t$, let us define
%
\begin{equation}
\label{defjt} J(t):=\inf \biggl\{s \Big| \int_0^s
u(v)\,\DD v\ge t \biggr\}.
\end{equation}
By construction, the process $(X_t)_{t\ge0}$ defined by
%
\begin{equation}
\label{defx} X_t:=A\bigl(J(t)\bigr)
\end{equation}
is a continuous time random walk on $\bbZ_+$ which jumps up and down
with rate $1/2$.
From the definition, we have
%
\begin{equation}
X_0=A(0):= N\mathcal H_0.
\end{equation}
Note that from Proposition~\ref{fluctouse}, and the definition of $\cH
_0$ we have
%
\begin{equation}
\label{fructose} \bbP \bigl[A(0)\ge2N^{3/2}\log\log N \bigr]\le(\log
N)^{-1}.
\end{equation}

To estimate $\tau$, we have to control the evolution of
$X_t$ (using standard properties of the random walk) and that of $u(t)$
(using the properties of proved in Section~\ref{repair}).

\subsection{\texorpdfstring{Multiscale analysis.}{Multiscale analysis}}\label{craak}

To have the best possible control on $u(t)$, we need to perform a
multi-scale analysis.
We construct a sequence of intermediate stopping time $(\tau_i)_{i\ge
0}$ as follows:
%
\begin{equation}
\tau_i:=\inf \bigl\{t\ge0 \mid A(t)\le N^{3/2}2^{-i}
\bigr\}.
\end{equation}
We set $\tau_{-1}:=0$ for convenience.
We are interested in $\tau_i$ for $i\in\{0,\dots,\break \lceil(\log_2
N)/ 2 \rceil\}$ where $\log_2(\cdot):=\log(\cdot)/\log(2)$ denotes
the logarithm in base $2$.
To bound the value of $\tau$, we bound the value of each
$\gD\tau_i=\tau_{i}-\tau_{i-1}$ for $i\le\lceil(\log_2 N)/2
\rceil$ and that of $\tau-\tau_{\lceil(\log_2 N)/2 \rceil}$.
The way to do this is:
\begin{longlist}[(ii)]
\item[(i)] First, we prove a bound for the analog of the $\Delta\tau
_i$ for the process
$X_t$ defined in \eqref{defx}.
\item[(ii)] Second, we prove a bound for $u(t)$ which is valid in the
interval $[\tau_{i-1},\tau_i)$.
\end{longlist}
For step (i), let us define
%
\begin{eqnarray}
\label{freddo} \cT_i&:=&\int_{\tau_{i-1}}^{\tau_i}
u(t) \,\DD t,
\nonumber\\[-8pt]\\[-8pt]\nonumber
\cT_\infty&:=&\int_{\tau_{\lceil(\log_2 N)/2 \rceil}}^{\tau} u(t) \,\DD t.
\end{eqnarray}
It follows from standard properties of the random walk and from \eqref
{fructose} that
we have the following.

\begin{lemma}\label{cromican}
We have the following estimates:
%
\begin{eqnarray}
\qquad \bbP \bigl[ \exists i\in\bigl\{0,\dots, \bigl\lceil(\log_2 N)/2
\bigr\rceil\bigr\}, \cT_i \ge3^{-i} N^3 (\log
N)^{1/4} \bigr] &\le& (\log N)^{-1/10},
\nonumber\\[-8pt]\\[-8pt]\nonumber
\bbP \bigl[\cT_\infty\ge N^2 (\log N)^{1/4}
\bigr]&\le&(\log N)^{-1/10}.
\end{eqnarray}
\end{lemma}

For more details, we refer to the proof of \cite{cfLac2}, Lemma 6.1.

Step (ii) is more delicate, because we cannot get a good bound on $u$
which is uniform in time.
For instance, we need to prove that most of the time $u(t)$ is of order
$N$ but we know that
just before $\tau$ we have $u(t)=4$. Hence, we will prove a different
bound for each value of $i$.
The bound is valid most of the time, and we will need to check that the
small fraction of time during which it does not hold can be dealt with
in the computations.
Recalling \eqref{defh}, we set
%
\begin{equation}
\mathcal H(t):= \max \bigl( H_{\alpha,t}\bigl(\xi^1_t
\bigr)+H_{\alpha,t}\bigl(\xi ^2_t\bigr), \sqrt{N}
\bigr).
\end{equation}
We notice that from the definition
%
\begin{equation}
\max_{x\in\bbZ} \bigl(\xi^2_t(x)-
\xi^1_t(x) \bigr)\le\mathcal \cH(t).
\end{equation}
Using this information, we can get the following control on $u$ [recall
\eqref{defce}]:

\begin{lemma}\label{fromzea}
If $\xi^1_t\in\cE$, we have
%
\begin{equation}
\label{grotsd} u(t)\ge\frac{1} 8 \min \biggl( N, \frac{A(t)}{\mathcal\cH(t)}
\biggr).
\end{equation}
\end{lemma}

The\vspace*{1pt} proof is identical to the one of \cite{cfLac2}, Lemma 6.3. Note
that thanks to Proposition \eqref{iamzea} and our assumption $\alpha
\le2 N^{-3/7}$,
the inequality \eqref{grotsd} is valid up to time $N^3$ (which is much
more than what we need) with high probability.
To make this bound on $u$ useful, we need to show that most of the time
$\cH(t)$ is not too large.

\begin{lemma}\label{coniduam}
For any $T\ge0$,
%
\begin{equation}
\bbP \biggl[ \int^T_0 \ind_{\{ \cH(t)\ge\sqrt{N} \log\log N \}
}\,\DD
t \ge T (\log N)^{-4} \biggr] \le(\log N)^{-1}.
\end{equation}
\end{lemma}

\begin{pf}
It follows from \eqref{fluctouse} that for $N$ sufficiently large, for
any $t\ge0$
%
\begin{equation}
\bbP \bigl[ \cH(t)\ge\sqrt{N} \log\log N \bigr]\le(\log N)^{-5}.
\end{equation}
Then the result follows by using the Markov property for the integrated
inequality.
\end{pf}

\begin{pf*}{Proof of Proposition~\ref{nahnou}}
Set
%
\begin{eqnarray}
\mathcal A&:=& \bigl\{\forall t\le N^3, \xi^1_t
\in\cE \bigr\},\nonumber
\\
\mathcal B&:=& \biggl\{\int^{T}_0
\ind_{\{ \cH(t)\ge\sqrt{N}\log
\log N\}}\,\DD t\le T(\log N)^{-5} \biggr\},
\\
\mathcal C&:=& \bigl\{\cT_i\le3^{-i} N^3 (\log
N)^{1/4} \bigr\} \cap \bigl\{\cT_\infty\le N^2(\log
N)^{1/4} \bigr\},\nonumber
\end{eqnarray}
where
%
\begin{eqnarray}
T&:=& N^2\sqrt{\log N}.
\end{eqnarray}
We assume also that $N$ is large enough so that from Proposition~\ref
{iamzea} and Lemmas~\ref{cromican} and~\ref{coniduam} we have
%
\begin{equation}
\bbP[\cA\cap\cB\cap\cC] \ge1-(2\log\log N)^{-1}.
\end{equation}
Hence, the results follows if we can prove that
%
\begin{equation}
\label{letrucs} \{\cA\cap\cB\cap\cC\}\subset\{\tau\le T\}.
\end{equation}
We split the proof of \eqref{letrucs} in two statements.
We want to show first that on the event $\cA\cap\cB\cap\cC$
%
\begin{equation}
\label{ige1} \tau-\tau_{\lceil\log_2 N/2 \rceil}\le(\log N)^{1/4}
N^2,
\end{equation}
and then that
%
\begin{equation}
\label{ige2} \forall i\in\bigl\{0,\dots, \bigl\lceil(\log_2 N)/2
\bigr\rceil\bigr\}, \qquad(\tau _{i}-\tau_{i-1})
\le(i+1)^{-2} N^2(\log N)^{1/3}.\hspace*{-30pt}
\end{equation}
These inequalities combined give
%
\begin{equation}
\qquad\tau\le(\log N)^{1/4} N^2 +\sum
_{i=0}^K (i+1)^{-2}N^2(\log
N)^{1/3}\le N^2 \sqrt{\log N}.
\end{equation}
Note that \eqref{ige1} is an immediate consequence of $\cC$ as
%
\begin{equation}
\cT_\infty=\int_{\tau_K}^\tau u(t)\,\DD t\ge
\tau-\tau_K.
\end{equation}
Let us turn to \eqref{ige2}.
Let us assume that the statement is false and set
%
\begin{eqnarray}
i_0 &:=& \min \bigl\{i \in\bigl\{0,\dots,
\nonumber\\[-8pt]\\[-8pt]\nonumber
&&{}  \bigl\lceil(\log_2
N)/2 \bigr\rceil\bigr\} \mid (\tau_{i}-\tau_{i-1})>
(i+1)^{-2}N^2(\log N)^{1/3} \bigr\}.
\end{eqnarray}
The definition of $i_0$ implies that
%
\begin{equation}
\label{bet} \tau_{i_0-1}+ (i_0+1)^{-2}N^2(
\log N)^{1/3} \le T.
\end{equation}
From $\cB$, we have [using \eqref{bet} to obtain the second inequality]
%
\begin{eqnarray}
\label{cardrive}
&& \int_{\tau_{i_0-1}}^{\tau_{i_0}} \ind_{\{H(t)\le\sqrt{N} \log
\log N \}}\nonumber
\\
&&\qquad \ge\int_{\tau_{i_0-1}}^{\tau_{i_0-1}+ (i_0+1)^{-2} N^2(\log
N)^{1/3}} \ind_{\{\cH(t)\le\sqrt{N} \log\log N \}} \,\DD t\nonumber
\\
&&\qquad =(i_0+1)^{-2} N^2(\log N)^{1/3}
\nonumber\\[-8pt]\\[-8pt]
&&\quad\qquad{}- \int
_{\tau_{i_0-1}}^{\tau_{i_0-1}+
(i_0+1)^{-2} N^2(\log N)^{1/3}} \ind_{\{\cH(t)> \sqrt{N} \log\log
N \}} \,\DD t\nonumber
\\
&&\qquad \ge(i_0+1)^{-2} N^2(\log N)^{1/3} -
N^2 (\log N)^{-3}\nonumber
\\
&&\qquad \ge\frac
{1}{2}(i_0+1)^{-2}
N^2(\log N)^{1/3}.\nonumber
\end{eqnarray}
For all $t\le\tau_{i_0}$, we have $A(t)\ge N^{3/2} 2^{-i_0}$, and
thus using Lemma~\ref{fromzea} and the assumption that $\cA$ holds,
%
\begin{eqnarray}
u(t) &\ge& \frac{1} 8 \min \biggl(N, \frac{A(t)}{\max(\cH
(t),N^{1/2})} \biggr)
\nonumber\\[-8pt]\\[-8pt]\nonumber
&\ge&
\frac{N^{3/2} 2^{-i_0}}{8\sqrt{N} \log\log
N}\ind_{\{\cH(t)\le\sqrt{N} \log\log N \}}.
\end{eqnarray}
From \eqref{cardrive},
%
\begin{eqnarray}
\cT_{i_0}&=&\int_{\tau_{i_0-1}}^{\tau_{i_0}} u(t)\,\DD t \ge
\frac{N
2^{-i_0} }{8 (\log\log N)} \int_{\tau_{i_0-1}}^{\tau_{i_0}}
\ind_{\{\cH(t)\le\sqrt{N} \log\log N \}} \,\DD t
\nonumber\\[-8pt]\\[-8pt]\nonumber
&\ge&(i_0+1)^{-2} 2^{-i_0} \frac{N^3(\log N)^{1/3}}{16 \log\log
N}>3^{-i_0}N^3(
\log N)^{1/4}.
\end{eqnarray}
This brings a contradiction to $\cC$ (if $N$ is large enough) and ends
the proof of \eqref{letrucs}.
\end{pf*}

\section{\texorpdfstring{Proof of Proposition \protect\ref{densitevol}.}{Proof of Proposition 4.3}} \label{relax}

To prove the result, we will try to control the derivative in $t$ of
the total variation distance that we have to bound.

Note that
$\llVert   \nu^{\alpha}_{t}-\nu^{\alpha e^{-\gl_{N}t}}\rrVert  _{\mathrm{TV}}$ is always
differentiable on the right.
This comes from the fact that for each $\eta\in\gO_N$, both $\nu
^{\alpha}_{t}(\eta)$ and $\nu^{\alpha e^{-\gl_Nt}}$ are
differentiable.
With a small abuse of notation, we use $\partial_t$ to denote the
right derivative.
Our method to prove Proposition~\ref{densitevol} relies on getting a
bound on the derivative valid for all $\alpha\le1$.
More precisely, we want to prove
%
\begin{equation}
\label{groo} \partial_t\bigl\llVert \nu^{\alpha}_{t}-
\nu^{\alpha e^{-\gl_Nt}}\bigr\rrVert _{\mathrm{TV}}\le C_1
\alpha^3 N^{-2} e^{-3\gl_N t}+ C_2
\alpha^2 N^{-3/2} e^{-2\gl_N t}.
\end{equation}
Indeed, once integrated this gives
%
\begin{equation}
\sup_{t\ge0}\bigl\llVert \nu^{\alpha}_{t}-
\nu^{\alpha e^{-\gl
_Nt}}\bigr\rrVert _{\mathrm{TV}}\le C_3 \bigl(
\alpha^3+\alpha^2 N^{1/2} \bigr),
\end{equation}
which is equivalent to our result.

Let us first perform a simple computation to show that it is sufficient
to prove \eqref{groo} in the case $t=0$.
Using the triangular inequality and the fact that the semi-group shrinks
the total-variation distance, we have for any positive $\gep$,
%
\begin{eqnarray}
&& \bigl\llVert \nu^{\alpha}_{t+\gep}-\nu^{\alpha e^{-\gl_N(t+\gep
)}}\bigr\rrVert
_{\mathrm{TV}}\nonumber
\\
&&\qquad \le \bigl\llVert \nu^{\alpha}_{t+\gep}-
\nu^{\alpha e^{-\gl_Nt}}_\gep \bigr\rrVert _{\mathrm{TV}}+ \bigl\llVert
\nu^{\alpha e^{-\gl_Nt}}_\gep-\nu^{\alpha e^{-\gl
_N(t+\gep)}} \bigr\rrVert
_{\mathrm{TV}}
\\
&&\qquad \le \bigl\llVert \nu^{\alpha}_{t}-\nu^{\alpha e^{-\gl_Nt}}\bigr
\rrVert _{\mathrm{TV}} + \bigl\llVert \nu^{\alpha e^{-\gl_Nt}}_\gep-
\nu^{\alpha e^{-\gl
_N(t+\gep)}} \bigr\rrVert _{\mathrm{TV}}.\nonumber
\end{eqnarray}

Hence,
%
\begin{equation}
\label{rocknroll} \partial_t\bigl\llVert \nu^{\alpha}_{t}-
\nu^{\alpha e^{-\gl_Nt}}\bigr\rrVert _{\mathrm{TV}} \le\partial_{\gep} \bigl
\llVert \nu^{\alpha e^{-\gl_Nt}}_\gep -\nu^{\alpha e^{-\gl_N(t+\gep)}} \bigr\rrVert
_{\mathrm{TV}}\mid _{\gep=0}.
\end{equation}

Note that the right-hand side is simply equal to
\[
\partial_s \bigl\llVert \nu^{\alpha'}_{s}-
\nu^{\alpha' e^{-\gl
_Ns}}\bigr\rrVert _{\mathrm{TV}} \mid _{s=0}
\]
for $\alpha'=\alpha e^{-\gl t}$.
Hence, to prove \eqref{groo} it is sufficient to show that for all
$\alpha\le1$
%
\begin{equation}
\label{daabound} \partial_{t} \bigl\llVert \nu^{\alpha}_t
-\nu^{\alpha e^{-\gl_Nt}} \bigr\rrVert _{\mathrm{TV}}\mid_{t=0}\le
C_1\alpha^3 N^{-2}+ C_2
\alpha^2 N^{-3/2}.
\end{equation}
We let $g^{\alpha}_t$ denote the density of $\nu^{\alpha}_{t}$, and
$g^\alpha$ that of $\nu^{\alpha}$.
Recall that we have
%
\begin{equation}
g^{\alpha}(\eta):=\frac{e^{\alpha a(\eta)}}{\mu_N(e^{\alpha a(\eta)})},
\end{equation}
where
%
\begin{equation}
a(\eta):=\sum_{x\in\bbZ_{2N}}\eta(x)\bsin(x).
\end{equation}
We have
%
\begin{equation}
\label{mikp} \partial_{t} \bigl\llVert \nu^{\alpha}_t
-\nu^{\alpha e^{-\gl_Nt}} \bigr\rrVert _{\mathrm{TV}} = \mu_N \bigl
\llvert \partial_t \bigl( g^{\alpha}_t(\eta)-
g^{\alpha e^{-\gl_Nt}}(\eta) \bigr)\mid _{t=0} \bigr\rrvert.
\end{equation}
We compute the derivatives of $g^{\alpha}_t$ and $g^{\alpha e^{-\gl
_Nt}}(\eta)$ separately.
We have
%
\begin{equation}
\label{crcri} \partial_{t} g^{\alpha e^{-\gl_Nt}}(\eta)
\mid_{t=0}= \alpha\gl_N g^\alpha(\eta) \bigl[- a(
\eta)+\nu^\alpha\bigl(a(\eta)\bigr) \bigr].
\end{equation}
The other term requires more work, and we have
%
\begin{equation}
\partial_{t}g^{\alpha}_{t}(\eta)
\mid_{t=0}=\mathcal L g= \sum_{x\in\bbZ_{2N}}
g^{\alpha}\bigl(\eta^x\bigr)-g^{\alpha}(\eta).
\end{equation}
Recall \eqref{defgrad}. We have
%
\begin{equation}
g^{\alpha}\bigl(\eta^x\bigr)-g^{\alpha}(
\eta)=g^{\alpha}(\eta) \bigl[\exp \bigl(-\alpha\nabla\bsin(x)\nabla\eta(x)
\bigr)-1 \bigr].
\end{equation}
Performing a Taylor expansion of the exponential, we have
%
\begin{equation}
\label{protos} \qquad\mathcal L g:= g^{\alpha}(\eta) \biggl[ - \alpha \biggl(
\sum_{x\in
\bbZ_{2N}}\nabla\bsin(x)\nabla\eta(x) \biggr)+
\frac{\alpha^2}{2}G(\eta,N)+R(\eta,N) \biggr],
\end{equation}
where $(\alpha^2/2)G(\eta,N)$ is the second term in the Taylor expansion
%
\begin{equation}
G(\eta,N):=\sum_{x\in\bbZ_{2N}}\bigl(\nabla\overline{\sin}
(x)\bigr)^2\bigl(\nabla\eta(x)\bigr)^2,
\end{equation}
and $R(\eta,N)$ is the Taylor rest
%
\begin{eqnarray}
R(\eta,N) &:=& \sum_{x\in\bbZ_{2N}} \biggl( e^{-\alpha\nabla\bsin
(x)\nabla\eta(x)}-1
+\alpha\nabla\overline{\sin} (x)\nabla\eta(x)
\nonumber\\[-8pt]\\[-8pt]\nonumber
&&{} -\frac{\alpha
^2}{2} \bigl(\nabla
\overline{\sin} (x)\bigr)^2\bigl(\nabla\eta(x)\bigr)^2
\biggr).
\end{eqnarray}
The first term in the RHS of \eqref{protos} can be simplified
using summation by part and the fact that $\overline{\sin}$ is an
eigenfunction of $\gD$. We have
%
\begin{equation}
\label{walsh} \sum_{x\in\bbZ_{2N}}\nabla\overline{\sin}(x)
\nabla\eta(x) =-\sum_{x\in\bbZ_{2N}}\gD\overline{\sin}(x)
\eta(x)=\gl_
{N}a(\eta).
\end{equation}
According to Taylor's formula, one has for all $\alpha<1$, for an
adequate choice of constant $C_1$
%
\begin{equation}
\label{grimos} \bigl\llvert R(\eta,N) \bigr\rrvert \le\frac{e^{2\alpha}\alpha^3}{6}\sum
_{x\in\bbZ
_{2N}} \bigl\llvert \nabla\overline{\sin}(x)\nabla
\eta(x)\bigr\rrvert ^3\le C_1\alpha^3
N^{-2},
\end{equation}
where in the last inequality we have used that $\llvert  \nabla\eta(x)\rrvert  \le2$
and that
%
\begin{equation}
\label{pisurn} \bigl\llvert \nabla\overline{\sin} (x)\bigr\rrvert =2\biggl\llvert
\sin \biggl(\frac{\pi}{2N} \biggr)\cos \biggl(\frac{\pi x}{N}+
\frac{\pi}{2N} \biggr)\biggr\rrvert \le\frac{\pi}{N}.
\end{equation}
Combining \eqref{mikp} with \eqref{crcri} and \eqref{walsh}, we obtain
%
\begin{eqnarray}
\partial_{t} \bigl\llVert \nu^{\alpha}_t -
\nu^{\alpha e^{-\gl_Nt}} \bigr\rrVert _{\mathrm{TV}}&\le& \nu_{\alpha}\biggl
\llvert \frac{\alpha^2}{2}G(\eta,N) +R(\eta,N)-\alpha\gl_N
\nu^\alpha\bigl(a(\eta)\bigr)\biggr\rrvert\nonumber
\\
&\le& \nu_{\alpha} \biggl\llvert R(\eta,N) -\alpha\gl_N
\nu^\alpha\bigl(a(\eta)\bigr)+\frac{\alpha^2}{2}\nu_{\alpha}
\bigl(G(\eta,N)\bigr) \biggr\rrvert
\\
&&{}+ \frac{\alpha^2}{2} \nu^{\alpha}\bigl\llvert G(\eta,N)-
\nu^{\alpha} \bigl( G(\eta,N) \bigr)\bigr\rrvert.\nonumber
\end{eqnarray}
To estimate the first term, we note that as
%
\begin{equation}
\mu_N \bigl( \partial_t \bigl( g^{\alpha}_t(
\eta)- g^{\alpha e^{-\gl_Nt}}(\eta) \bigr)\mid _{t=0} \bigr)=0,
\end{equation}
we have
%
\begin{equation}
\nu_{\alpha} \biggl( R(\eta,N) -\alpha\gl_N
\nu^\alpha\bigl(a(\eta)\bigr)+\frac{\alpha^2}{2}G(\eta,N) \biggr)=0.
\end{equation}
Hence, from \eqref{grimos}
%
\begin{eqnarray}
&& \nu_{\alpha} \biggl\llvert R(\eta,N) -\alpha\gl_N
\nu^\alpha\bigl(a(\eta)\bigr)+\frac{\alpha^2}{2}\nu_{\alpha}
\bigl(G(\eta,N)\bigr) \biggr\rrvert
\nonumber\\[-8pt]\\[-8pt]\nonumber
&&\qquad = \nu_{\alpha} \bigl\llvert R(\eta,N) -\nu^{\alpha}\bigl(R(\eta,N)
\bigr)\bigr\rrvert \le C_1 \alpha^3 N^{-2}.
\end{eqnarray}
To estimate the second term, we use Proposition~\ref{lipitz}.
The reader can check that the Lipshitz norm of $G$ [cf. \eqref
{lipdef}] of $G$ satisfies
%
\begin{equation}
\bigl\llVert G(\cdot,N)\bigr\rrVert _{\lip}\le8 \pi^2
N^{-2}
\end{equation}
and hence that for an adequate choice of $C_2>0$
%
\begin{equation}
\nu^{\alpha}\bigl\llvert G(\eta,N)- \nu^{\alpha} \bigl( G(\eta,N)
\bigr)\bigr\rrvert \le C_2 N^{-3/2}.
\end{equation}
This completes the proof of \eqref{daabound}.

\begin{appendix}
\section{\texorpdfstring{Proof of technical statements on $\nu^\alpha$.}{Proof of technical statements on nualpha}}\label{techos}

\subsection{\texorpdfstring{Proof of Proposition \protect\ref{rightdensity}.}{Proof of Proposition 6.1}}

Note that if $\mu_N$ was replaced by the uniform measure on $\{-1,1\}
^{\bbZ_{2N}}$ (without the constraint of having $N$ particles)
then $\nu^\alpha$ would be a product of independent Bernoulli, and
the statement would be trivial to prove.

What we have to control is that the constraint on the number of
particles does not affect the mean too much. To do so, we perform an
expansion of the partition function
according to the value of $\eta(x)$ to show that the ratio of the
partition function restricted to the event $\eta(x)=+1$ and $\eta
(x)=-1$, respectively, is close to\vadjust{\goodbreak}
$\exp(2\alpha\bsin(x))$. To this purpose, we introduce the quantity
%
\begin{equation}
Z(x):=\frac{\mu_N (e^{\alpha\sum_{y\in\bbZ_{2N}\setminus\{
x\}}\eta(y) \bsin(y)}  \mid \eta(x)=+1 )}{
\mu_N (e^{\alpha\sum_{y\in\bbZ_{2N}\setminus\{x\}}\eta(y)
\bsin(y)}  \mid \eta(x)=-1 )}.
\end{equation}
We have
\begin{eqnarray}
\label{croqui}
&& \nu^{\alpha}\bigl(\eta(x)\bigr)\nonumber
\\
&&\qquad = \nu^{\alpha}\bigl(
\eta(x)=+1\bigr)-\nu^{\alpha
}\bigl(\eta(x)=-1\bigr)
\nonumber\\[-8pt]\\[-8pt]\nonumber
&&\qquad = \frac{\mu_N (e^{\alpha\sum_{y\in\bbZ_{2N}}\eta(y) \bsin
(y)}  \llvert    \eta(x)=+1 )-\mu_N (e^{\alpha\sum_{y\in\bbZ
_{2N}}\eta(y) \bsin(y)}  \rrvert    \eta(x)=-1 )}{
\mu_N (e^{\alpha\sum_{y\in\bbZ_{2N}}\eta(y) \bsin(y)}  \llvert
 \eta(x)=+1 )+\mu_N (e^{\alpha\sum_{y\in\bbZ
_{2N}}\eta(y) \bsin(y)}  \rrvert    \eta(x)=-1 )}\nonumber
\\
&&\qquad = \frac{e^{2\alpha\bsin(x)}Z(x)-1}{e^{2\alpha\sin(x)}Z(x)+1}.\nonumber
\end{eqnarray}
Hence, what we must check to prove \eqref{dsss} is that $Z(x)$ is very
close to one.
Now note that we can obtain a coupling of $\mu_N ( \cdot \mid\eta(x)=-1 )$ and $\mu_N ( \cdot \mid  \eta(x)=+1
)$ in the following manner:
take $\eta^1$ with distribution $\mu_N ( \cdot \mid  \eta
(x)=-1 )$, choose $y$ uniformly at random (and independent of
$\eta^1$ in $\{z  \mid  \eta^1(z)=+1 \}$ and let
$\eta^2$ be obtained from $\eta^1$ by exchanging the value at $x$ and
$y$ (which are $+1$ and $-1$, resp.).
A~consequence of this coupling is that
%
\begin{eqnarray}
&& Z(x)
\nonumber\\[-6pt]\\[-14pt]\nonumber
&&\!\qquad :=\frac{ \sklfrac{1}{N} \mu_N (\sum_{y\in\{z  \mid   \eta
(z)=+1 \}} e^{\alpha\sum_{w\in\bbZ_{2N}\setminus\{x\}}\eta(w)
\bsin(w)-2\alpha\bsin(y)}  \mid   \eta(x)=-1 ) }{ \mu_N (e^{\alpha\sum_{w\in\bbZ_{2N}\setminus\{x\}}\eta(w) \bsin(w)}  \mid  \eta(x)=-1 )},\hspace*{-10pt}
\end{eqnarray}
and hence we can deduce from it
%
\begin{eqnarray}
\label{five} Z(x) &=& \nu^{\alpha} \biggl( \frac{1}{N}\sum
_{y\in\bbZ_{2N}}\frac{1+\eta
(y)}{2}e^{-2\alpha\bsin(y)} \Big| \eta(x)=-1
\biggr)
\nonumber\\[-8pt]\\[-8pt]\nonumber
&=& 1+ \nu^{\alpha} \biggl( \frac{1}{N}\sum
_{y\in\bbZ_{2N}}\frac{1+\eta
(y)}{2} \bigl(e^{-2\alpha\bsin(y)}-1 \bigr) \Big|
\eta(x)=-1 \biggr).
\end{eqnarray}

Note that with this expression it is not hard to check that
$\llvert  Z(x)-1\rrvert  \le e^{2\alpha}-1$.
However, to get a sharper estimate, we must have a good control on
$ \nu^{\alpha} (\eta(y) \mid \eta(x)=-1  )$.
We obtain it by pushing the expansion one step further. We set
%
\begin{equation}
Z'(x,y)=\frac{\mu_N (e^{\alpha\sum_{z\in\bbZ_{2N}\setminus\{
x,y\}} \bsin(z)}  \mid  \eta(x)=+1,\eta(y)=-1 )}{
\mu_N (e^{\alpha\sum_{z\in\bbZ_{2N}\setminus\{x,y\}} \bsin
(z)}  \mid  \eta(x)=-1,\eta(y)=-1 )}.
\end{equation}
Similar to \eqref{croqui}, we obtain that
%
\begin{equation}
\label{croqui2} \nu^{\alpha}\bigl(\eta(y) \mid \eta(x)=-1\bigr)=
\frac{N}{N-1}\frac
{e^{2\alpha\bsin(y)}Z'(x,y)-1}{e^{2\alpha\bsin(y)}Z'(x,y)+1}.
\end{equation}
Like for \eqref{five}, we have an alternative expression for $Z'$
%
\begin{eqnarray}
&& Z'(x,y)
\nonumber\\[-6pt]\\[-14pt]\nonumber
&&\qquad =1+\nu^{\alpha} \biggl( \frac{1}{N}\sum
_{z\in\bbZ
_{2N}}\frac{1+\eta(z)}{2} \bigl(e^{-2\alpha\bsin(z)}-1
\bigr) \Big|  \eta(x)=-1, \eta(y)=-1 \biggr).\hspace*{-20pt}
\end{eqnarray}
Hence, we have
%
\begin{equation}
\bigl\llvert Z'(x,y)-1\bigr\rrvert \le e^{2\alpha}-1,
\end{equation}
and from \eqref{croqui2}, we deduce that for some positive constant $C_1$
%
\begin{equation}
\bigl\llvert \nu^{\alpha}\bigl(\eta(y) \mid \eta(x)=-1\bigr) \bigr\rrvert
\le C_1 \biggl( \frac{1}{N}+\alpha \biggr).
\end{equation}
Hence, we have
%
\begin{eqnarray}
\label{grimm}
\qquad\bigl\llvert Z(x)-1\bigr\rrvert &\le&\frac{1}{2N}\biggl\llvert
\sum_{y\in\bbZ_{2N}\setminus\{x\}
} \bigl(e^{-2\alpha\bsin(y)}-1 \bigr)\biggr
\rrvert
\nonumber\\[-8pt]\\[-8pt]\nonumber
&&{}+\frac{1}{2N}\sum_{y\in\bbZ_{2N}\setminus\{x\}} \bigl\llvert \nu
^{\alpha}\bigl(\eta(y) \mid \eta(x)=-1\bigr) \bigl(e^{-2\alpha\bsin(y)}-1
\bigr)\bigr\rrvert.
\end{eqnarray}
Performing a Taylor expansion up to the second order in $\alpha$ we
obtain (recall $\alpha\le1$)
%
\begin{equation}
\biggl\llvert \sum_{y\in\bbZ_{2N}\setminus\{x\}} \bigl(e^{-2\alpha\bsin
(y)}-1
\bigr) \biggr\rrvert \le2\alpha\bigl\llvert \bsin(x) \bigr\rrvert +
\frac{e N\alpha^2}{2}.
\end{equation}
The second term in the RHS of \eqref{grimm} can be bounded by
%
\begin{equation}
C_1 \biggl( \frac{1}{N}+\alpha \biggr) \bigl(e^\alpha-1
\bigr).
\end{equation}
Hence, we obtain
%
\begin{equation}
\bigl\llvert Z(x)-1\bigr\rrvert \le C_2 \bigl( \alpha^2
+ N^{-2}\bigr).
\end{equation}
And then the result can easily be deduced from \eqref{croqui}.

\subsection{\texorpdfstring{Proof of Proposition \protect\ref{fluctouse}.}{Proof of Proposition 6.2}}

The result follows from the combination of Proposition~\ref{sinusoid2}
which controls the fluctuation around the mean value $u^{\eta_0}(x,t)$
given an initial condition $\eta_0$
and the following statement, that the mean itself $u^{\eta_0}(x,t)$
does not fluctuate too much if $\eta_0$ has distribution
$\nu^\alpha$.

\begin{lemma}
There exists a constant $c$ such that for all $N>0$, for all $\llvert   \alpha \rrvert  \le N^{-1/4}$, and $t\ge0$
%
\begin{equation}
\label{fluctumean} \nu^{\alpha} \Biggl[ \exists x, y \in\bbZ_{2N},
\Biggl\llvert \sum_{z=x+1}^y
\bigl(u^\eta(x,t)-\alpha e^{-\gl_N t}\bsin (z) \bigr) \Biggr\rrvert
\ge s \sqrt{N} \Biggr]\le2e^{-cs^2}.\hspace*{-30pt}
\end{equation}
\end{lemma}

\begin{pf}
It is in fact sufficient to prove \eqref{fluctumean} for $t=0$, because
%
\begin{equation}
\label{contract} \qquad\max_{x,y}\sum_{z=x+1}^y
\bigl(u^{\eta}(z,t)-\alpha e^{-\gl_N
t}\bsin(z) \bigr)\le \max
_{x,y}\sum_{z=x+1}^y
\bigl(\eta_0(z)-\alpha\bsin(z)\bigr).
\end{equation}
Indeed, if one sets $v(x,t)$ to be the solution of the discrete-heat
equation on $\bbZ_{2N}$
with initial condition
\[
v_0(x):=\sum_{z=1}^x
\eta(z)-\alpha\bsin(z),
\]
then \eqref{contract} can be reformulated as
%
\begin{equation}
\max_{x,y} \bigl[v(t,y)-v(t,x) \bigr]\le\max
_{x,y} \bigl[ v_0(y)-v_0(x) \bigr]
\end{equation}
which is obviously true by contractivity of the heat equation.
Note that at the cost of losing a factor in the constant $c$, we can
restrict ourselves to proving that
%
\begin{equation}
\nu^\alpha \Biggl[ \exists y \in\bbZ_{2N}, \Biggl\llvert \sum
_{z=1}^y \bigl[\eta(z)- \nu^{\alpha}
\bigl(\eta(z)\bigr) \bigr] \Biggr\rrvert \ge4s \sqrt{N} \Biggr]
\le2e^{-cs^2}.
\end{equation}
We have used Proposition~\ref{rightdensity} and the assumption on
$\alpha$ to replace $\alpha\bsin(z)$ by $\nu^{\alpha}(\eta(z))$.
Let us introduce notation for the sum
%
\begin{equation}
S_{x,y}:=\sum_{z=1}^y \bigl(
\eta(z)- \nu^{\alpha}\bigl(\eta(z)\bigr) \bigr).
\end{equation}
We also set $p:= \lfloor\log_2 N\rfloor+1$.
For $s>0$, we set
%
\begin{eqnarray}
\mathcal J(s) &:=& \bigl\{\exists q\in\{1,\dots,p\}, \exists y\in\bigl\{1,\dots,
\bigl\lfloor 2N2^{-q}\bigr\rfloor\bigr\},
\nonumber\\[-8pt]\\[-8pt]\nonumber
&&{}\llvert S_{2^{q}(y-1), 2^{q}y}\rrvert \ge \bigl(\tfrac{3}{4}
\bigr)^{p-q} s\sqrt {N} \bigr\}.
\end{eqnarray}
By a simple dichotomy argument (see the proof of Proposition 4.1 in
\cite{cfLac2}), we have
%
\begin{equation}
\Biggl\{ \exists y \in\bbZ_{2N}, \Biggl\llvert \sum
_{z=1}^y \bigl(\eta(z) - 1/2 - \alpha\bsin(z)
\bigr) \Biggr\rrvert \ge4s \sqrt{N} \Biggr\}\subset\mathcal J(s).
\end{equation}
For $y$ and $p$ fixed,
$S_{2^{q}(y-1), 2^{q}y}(\eta)$ is a function which depends on $2^q$
coordinates and whose Lipshitz norm is smaller than
$2$. Hence, by Proposition~\ref{lipitz}, we have
%
\begin{eqnarray}
\nu^{\alpha} \bigl( \llvert S_{2^{q}(y-1), 2^{q}y}\rrvert \ge \bigl(
\tfrac
{3}{4} \bigr)^{p-q} s\sqrt{N} \bigr) &\le&2\exp
\bigl(-C_1 \bigl(\tfrac{9}{16} \bigr)^{p-q}2^{-q}s^2
N \bigr)
\nonumber\\[-8pt]\\[-8pt]\nonumber
&\le&2\exp \bigl(-2C_1 \bigl(\tfrac{9}{8}
\bigr)^{p-q}s^2 \bigr).
\end{eqnarray}

Hence, by a union bound, for an appropriate choice of
constant $C_2$ and for all $s>0$, we have
%
\begin{equation}
\nu^{\alpha} \bigl(\mathcal H(s) \bigr)\le2 \sum
_{q=1}^p 2^{p-q} \exp
\biggl(-2C_1 \biggl(\frac{9}{8} \biggr)^{p-q}s^2
\biggr) \le2\exp \bigl(-C_2 s^2 \bigr).
\end{equation}\upqed
\end{pf}

\subsection{\texorpdfstring{Proof of Proposition \protect\ref{iamzea}.}{Proof of Proposition 6.3}}

Set
%
\begin{equation}
\bar\cE:=\bigl\{\eta\in\gO_N \mid \exists(x,y) \in
\bbZ^2_{2N}, \# [x,y]\ge N^{1/4}\Rightarrow j(x,y,
\eta) \ge\tfrac{1}{3}\#[x,y] \bigr\}.\hspace*{-30pt}
\end{equation}
First, we notice that from the proof of \cite{cfLac1}, Lemma 6.2,
there exists a constant $C_1>0$ such that
%
\begin{equation}
\mu_N(\bar\cE)\le e^{ -C_1N^{1/4}}.
\end{equation}
Recall that $\nu^\alpha_t$ the law of $\eta_t$ starting from
distribution $\nu^\alpha$.
We have by the  Cauchy--Schwarz inequality
%
\begin{equation}
\bigl(\nu^\alpha_t(\bar\cE) \bigr)^2 \le
\mu_N(\bar\cE) \mu _N \biggl[ \biggl( \frac{\DD\nu^\alpha_t}{\DD\mu_N }
\biggr)^2 \biggr].
\end{equation}
Note\vspace*{1pt} that the term $\mu_N [  ( \frac{\DD\nu^\alpha
_t}{\DD\mu_N }  )^2  ]$ is decreasing in $t$, because the
semi-group of the Markov chain contracts the $l_2$ norm.
For $t=0$, we have
%
\begin{equation}
\mu_N \biggl[ \biggl( \frac{\DD\nu^\alpha}{\DD\mu_N } \biggr)^2
\biggr] \le\mu_N \bigl(e^{2\alpha a(\eta)} \bigr).
\end{equation}
Using Proposition~\ref{lipitz} to have Gaussian concentration for
$a(\theta)$,
we have for $N$ sufficiently large:
%
\begin{equation}
\mu_N \biggl[ \biggl( \frac{\DD\nu^\alpha}{\DD\mu_N } \biggr)^2
\biggr]\le\exp \bigl(100N \alpha^2 \bigr).
\end{equation}
Hence, we can conclude that there exists constant $C_2$ and $C_3$ such
that if $\alpha< C_2 N^{-3/8}$ for any $t$ we have
%
\begin{equation}
\label{cromic} \nu^\alpha_t(\bar\cE)= \bbP^{\nu^{\alpha}}[
\eta_t \in\bar\cE ]\le e^{-C_3N^{1/4}}.
\end{equation}
Now we have to move from this result to a result for all $t\le N^3$.
Note that starting from $\eta\notin\bar\cE$, one needs at least
$\frac{1}{12} N^{1/4}$ transitions in order to
jump out of $\cE$.
Hence, using union bound
%
\begin{eqnarray}
&& \bbP^{\nu^{\alpha}}\bigl[ \exists t\le N^3, \eta_t
\notin\cE\bigr]\nonumber
\\
&&\qquad \le \sum_{i=0}^{N^5}
\bbP^{\nu^{\alpha}}[\eta_{i/N^2} \in\bar\cE]
+ \sum_{i=1}^{N^5} \bbP \biggl[ (
\eta_t)_{t\in
[(i-1)/N^2,i/N^2]}
\\
&&\quad\qquad{} \mbox{performs more than }
\frac{1}{12} N^{1/4} \mbox{ transitions} \biggr].\nonumber
\end{eqnarray}
The first term is smaller than $e^{-C_3N^{1/4}}$; cf. \eqref{cromic}.
As for the second one, it is not difficult to check that the rate at
which transitions occur in the chain is bounded by $2N$,
and thus that for any $i$
%
\begin{equation}
\bbP \bigl[ (\eta_t)_{t\in[(i-1)/N^2,i/N^2]}\mbox{ performs more than
} \tfrac{1}{12} N^{1/4}\mbox{ transitions} \bigr]\le
e^{-N},\hspace*{-30pt}
\end{equation}
provided $N$ is large enough.

\section{\texorpdfstring{Concentration for Lipschitz function of particle systems.}{Concentration for Lipschitz function of particle systems}}

Given $f: \{0,1\}^\bbZ_{2N} \to\bbR$,
one sets $\llVert  f\rrVert  _{\lip}$ to be the Lipschitz norm of $f$ for the
Hamming distance
%
\begin{equation}
\label{lipdef} \llVert f\rrVert _{\lip}:=\max_{\eta, \eta' \in\{-1,1\}^{\bbZ_{2N}}  }
\frac{\llvert  f(\eta)-f(\eta')\rrvert  }{\sum_{x\in\bbZ_{2N}} \ind_{\{ \eta
(x)\ne\eta'(x)\}}}.
\end{equation}

\begin{proposition}\label{lipitz}
For any $f \{-1,1\}^{\bbZ_{2N}} \to\bbR$ we have
%
\begin{equation}
\mu_N \bigl( \bigl\llvert f- \mu_N(f)\bigr\rrvert \ge s
\bigr) \le2 \exp \biggl( -\frac
{s^2}{8(2N-1)\llVert   f \rrVert  ^2_{\lip}} \biggr).
\end{equation}
If the function $f$ only depends on $(\eta_x)_{x\in A}$ where $A$ is
fixed a subset of $\bbZ_{2N}$ of cardinal $k$ we have
%
\begin{equation}
\mu_N \bigl( \bigl\llvert f- \mu_N(f)\bigr\rrvert \ge s
\bigr) \le2 \exp \biggl( -\frac
{s^2}{8k\llVert   f \rrVert  ^2_{\lip}} \biggr).
\end{equation}

The result remains valid if $\mu_N$ is replaced by
a measure $\nu$ whose density with respect to $\mu_N$ is of the form
%
\begin{equation}
\frac{\DD\nu}{\DD\mu_N}:= \frac{e^{\sum_{x\in\bbZ_{2N}}
g(x)\eta(x)}}{\mu_N(e^{\sum_{x\in\bbZ_{2N}} g(x)\eta(x)})},
\end{equation}
where $g$ is an arbitrary function on $\bbZ_{2N}$.
\end{proposition}

\begin{pf}
We can without loss of generality assume that $\llVert  f\rrVert  _{\lip}=1$.
Now, we introduce the martingale $(M_i)_{i=0}^{2N-1}$ defined by
%
\begin{equation}
M_i(\eta):= \nu \bigl( f(\eta) \mid \bigl(\eta(x)
\bigr)_{x=1}^i \bigr).
\end{equation}

We are going to check that the increments of $M$ are bounded, that is,
%
\begin{equation}
\label{boudoud} \forall i \in\{ 0,\dots, 2N-2\}, \qquad\llvert
M_{i+1}-M_i\rrvert \le2
\end{equation}
and the proposition is then simply a consequence of Azuma's
concentration inequality \cite{cfAzuma}.

To check \eqref{boudoud}, we need to show that for any realization
$(\eta(x))_{x=1}^{i}$
one can couple $\eta^1$ and $\eta^2$ with law
%
\begin{eqnarray}
\nu_1&:=&\nu \bigl( \cdot \mid \bigl(\eta(x)\bigr)_{x=1}^{i},
\eta(i+1)=1 \bigr),
\nonumber\\[-8pt]\\[-8pt]\nonumber
\nu_2&:=&\nu \bigl( \cdot \mid \bigl(\eta(x)\bigr)_{x=1}^{i},
\eta(i+1)=-1 \bigr)
\end{eqnarray}
in\vspace*{1pt} a way that $(\eta^1-\eta^2)(x)$ has only two discrepancies, one at
$i+1$ and another one in where $\eta^1(x)=1-\eta^2(x)=0$.

Note that $\nu_1$ and $\nu_2$ can be considered as a measure on $\{
-1,1\}\to\{i+2,\dots,2N\}$, one which is concentrated on the
set of configurations with $k:=N-\sum_{x=1}^i \eta(x)-1$ particles,
and the other on the set of configuration with $k+1$ particles.
What one can do is to first draw $\eta^1$ according to $\nu_1$, and
then add a $1$ chosen at
random to the configuration to obtain $\eta^2$.
One $\eta^1$ is given, and we choose at random a site $X$ in $\{x\in
\{ i+2,\dots, 2N\} \mid \eta^1(x)=-1\}$ with distribution
%
\begin{equation}
\frac{e^{g(x)}}{\sum_{\{x\in\{ i+2,\dots,2N\}\mid\eta^1(x)=-1 \}} e^{g(x)}}.
\end{equation}
On can check that $\eta^2$ defined by
%
\begin{equation}
\eta^2(x):= \eta^1(x)+\ind_{\{X=x\}} -
\ind_{\{x=i+1\}},
\end{equation}
has distribution $\nu_2$.

For the case where $f$ depends only on $\eta_{\mid A}$, we can consider a
$k$-step martingale which unveils
at each step the state $\eta(x)$ of one $x\in A$.
\end{pf}
\end{appendix}

\section*{\texorpdfstring{Acknowledgments.}{Acknowledgements}}
The author is very much indebted to Milton Jara who suggested
to him the question of cutoff profile, and pointed out that the
relaxation of the first Fourier coefficient of $\eta$ should be
similar to an Ornstein--Uhlenbeck process
starting far from equilibrium.



\printaddresses
\end{document}